\documentclass[12pt]{article}

\usepackage{amsmath}
\usepackage{amssymb}
\usepackage{amsthm}
\usepackage{graphicx}

\newtheorem{theorem}{Theorem}[section]
\newtheorem{lemma}[theorem]{Lemma}

\theoremstyle{definition}
\newtheorem{definition}[theorem]{Definition}
\newtheorem{example}[theorem]{Example}

\newtheorem{proposition}[theorem]{Proposition}
\newtheorem{corollary}[theorem]{Corollary}

\newtheorem{algorithm}[theorem]{Algorithm}

\theoremstyle{remark}
\newtheorem{remark}[theorem]{Remark}

\numberwithin{equation}{section}

\def\F{\mathbb{F}}
\def\Z{\mathbb{Z}}

\def\GG{\mathcal{G}}

\def\deg{\,\mbox{deg}\,}

\def\vec#1{(#1_1,#1_2,\ldots,#1_n)}

\def\tangle{\rightrightarrows}

\def\r#1{\rho^{(#1)}}

\def\th{\theta}
\def\T{\mathcal{T}}
\def\m#1{\mathop{|}_{r_{#1}}}
\def\ou#1{\overline{\underline{#1}}}
\def\o#1{\stackrel{r_{#1}}\longrightarrow}

\newenvironment{centre}{\begin{center}}{\end{center}}

\makeatletter
\newfont{\footsc}{cmcsc10 at 8truept}
\newfont{\footbf}{cmbx10 at 8truept}
\newfont{\footrm}{cmr10 at 10truept}
\renewcommand{\ps@plain}{%
\renewcommand{\@oddfoot}{\footsc the electronic journal of combinatorics
  {\footbf 13} (2006), \#R71\hfil\footrm\thepage}}
\makeatother 

\title{Factorial Grothendieck Polynomials}

\author{Peter J. McNamara\\
\small Department of Mathematics and Statistics\\[-0.8ex]
\small University of Sydney, NSW 2006, Australia\\[-0.8ex]
\small \texttt{petermcn@maths.usyd.edu.au}}

\date{\small
Submitted: Aug 10, 2005;  Accepted: Jan 8, 2006; Published: Aug 10, 2006
}

\begin{document}
\maketitle

\begin{abstract}
In this paper, we study Grothendieck polynomials indexed by
Grassmannian permutations from a combinatorial viewpoint. We
introduce the factorial Grothendieck polynomials which are
analogues of the factorial Schur functions, study their
properties, and use them to produce a generalisation of a
Littlewood-Richardson rule for Grothendieck polynomials.
\end{abstract}

\section{Introduction}
Let $x=(x_1,\dots,x_n)$ be a set of variables, $\beta$ a parameter
and $\theta$ a skew Young diagram whose columns have at most $n$ boxes.
A set-valued $\theta$-tableau $T$ is obtained by
placing subsets of $[n]=\{1,\dots,n\}$ (a notation used throughout) into the boxes of $T$
in such a way that the rows weakly increase while the
columns strictly increase. More precisely,
in each cell $\alpha$ of $\theta$, place a non-empty set
$T(\alpha)\subset[n]$ so that if $\alpha$ is
immediately to the left of $\beta$ then $\max(T(\alpha))\leq
\min(T(\beta))$, while if $\alpha$ is immediately above $\beta$,
then $\max(T(\alpha))<\min(T(\beta))$.
An example of such a $(4,4,2,1)/(1)$-tableau is given by the
following:

\[
\setlength{\unitlength}{1.46em}
\begin{picture}(4,4)
\multiput(1,2)(0,1){3}{\line(1,0){3}}
\multiput(4,4)(-1,0){4}{\line(0,-1){2}}
\multiput(0,0)(1,0){2}{\line(0,1){3}} \put(2,2){\line(0,-1){1}}
\put(1,1){\line(1,0){1}} \multiput(0,0)(0,1){4}{\line(1,0){1}}
\put(1,1.28){678} \put(0.35,0.28){9} \put(0.35,2.28){1}
\put(3.35,2.28){8} \put(1.35,2.28){4} \put(3.35,3.28){7}
\put(0.175,1.28){26} \put(2.175,3.28){34} \put(1.175,3.28){23}
\put(2.175,2.28){57}
\end{picture}
\setlength{\unitlength}{1pt}
\]

Given a skew diagram $\th$,
the (ordinary) Grothendieck polynomial
$G_\th(x)$ is defined by
\begin{equation} \label{odefint}
G_\th(x)=\sum_T \beta^{|T|-|\theta|} \prod_{\substack{\alpha\in\th
\\ r\in T(\alpha)}} x_r
\end{equation} where the sum is over all set-valued $\th$-tableaux $T$.

In a different form, the Grothendieck polynomials were first
introduced by Lascoux and Sch\"utzenberger \cite{lassch} as
representatives for $K$-theory classes determined by structure
sheaves of Schubert varieties. Since then, their properties were
studied by Fomin and Kirillov \cite{fkmini,FK}, Lenart
\cite{Lenart} and Buch \cite{B:buch}. In particular, the latter paper
contains the above combinatorial description of Grothendieck
polynomials in terms of tableaux, similar to that for the Schur
polynomials. It is this formulation which we use as the basis for
our approach to the study of Grothendieck polynomials in this
paper.

The major focus of this paper is the introduction and study of
what we shall call the factorial Grothendieck polynomials. They
generalise (\ref{odefint}) by introducing a second set of
variables $(a_i)_{i\in\Z}$ and we define the factorial
Grothendieck polynomial in $n$ variables $\vec{x}$ by
\[
G_\th(x|a)=\sum_T \beta^{|T|-|\theta|}
\prod_{\substack{\alpha\in\th
\\ r\in T(\alpha)}} x_r\oplus a_{r+c(\alpha)}
\] where $c(\alpha)$ is the {\it content} of the cell $\alpha$,
defined by $c(i,j)=j-i,$ and again the sum is over all set-valued
$\th$-tableaux $T$.

These factorial Grothendieck polynomials specialise in two
different ways, firstly by setting $a_i=0$ for all $i$ to obtain
the ordinary Grothendieck polynomials, and secondly by setting
$\beta=0$ to obtain the factorial Schur polynomials as studied in
\cite{MS:main}. Of these two families of polynomials obtained via
specialisation, the theory and properties of the factorial
Grothendieck polynomials appear to mimic more closely that of the
factorial Schur polynomials.

It can be shown and indeed is shown in this paper (Theorem
\ref{basis}) that the factorial Grothendieck polynomials
$G_{\lambda}(x|a)$ with $\lambda$ running over the (non-skew)
partitions with length at most $n$ form a basis of the ring of
symmetric polynomials in $x_1,\dots,x_n$. Hence, we can define the
coefficients $c^\nu_{\th\mu}(\beta,a)$ by the expansion
\begin{equation}\label{deflrint}
G_\theta(x|a)G_\mu(x|a)=\sum_\nu
c^\nu_{\th\mu}(\beta,a)\,G_\nu(x|a).
\end{equation}

In order to obtain a rule describing these coefficients, we
closely follow the method of Molev and Sagan
\cite{MS:main},\footnote{We are grateful to Anatol Kirillov for
suggesting to apply this method to the Grothendieck polynomials.}
exploiting the similarities between the factorial Grothendieck
polynomials and factorial Schur polynomials. This approach relies
on properties peculiar to the factorial versions of the
polynomials which enable a recurrence relation for the
coefficients to be determined, though there are also some
characteristics unique to the Grothendieck case, most notably in
section 4.2.

We present three solutions to the recurrence relation obtained for
the coefficients. The first of these is a general formula where
$G_\theta(x|a)$ in (\ref{deflrint}) is replaced by an arbitrary
symmetric polynomial. The second is a full solution in the case
where $\theta$ has no two boxes in the same column, which is
essentially a Pieri rule for factorial Grothendieck polynomials.
The third solution is a partial rule for arbitrary $\theta$
obtained by specialising certain variables to zero.

Out of the third solution, an application of the theory of
factorial Grothendieck polynomials to that of ordinary
Grothendieck polynomials is obtained. This consists of a
combinatorial rule for the calculation of the coefficients
$c_{\theta\mu}^{\nu}(\beta,0)$, generalising a previous result of
Buch \cite{B:buch}. In order to formulate the rule, define the
{\it column word\/} of a set-valued tableau $T$ as the sequence
obtained by reading the entries of $T$ from top to bottom in
successive columns starting from the right most column with the
rule that the entries of a particular box are read in the
decreasing order. As an example, the column word of the tableau
depicted earlier in the introduction is $7843753248761629$.

We write $\lambda\to\mu$ if $\mu$ is obtained by adding one box to
$\lambda.$ If $r$ is the row number of the box added to $\lambda$
to create $\mu$ then write $\lambda\stackrel{r}\to\mu$.
A set-valued tableau $T$ {\it fits\/}
a sequence $R(\mu,\nu)$ of partitions
$$
\mu=\r{0}\o{1}\r{1}\o{2}\cdots\o{l}\r{l}=\nu
$$
if the column word of $T$
coincides with $r_1\dots r_l$. With this notation, we have

\bigskip
\noindent {\bf Theorem}\quad {\it The coefficient
$c_{\theta\mu}^\nu(\beta,0)$ is equal to
$\beta^{|\nu|-|\mu|-|\theta|}$ times the number of set-valued
$\theta$-tableaux $T$ such that $T$ fits a sequence $R(\mu,\nu)$.
}

\bigskip

In the particular cases where $\theta=\lambda$ is normal, or
$\mu=\emptyset$, our rule coincides with the one previously given
by Buch \cite{B:buch}. Note also that if $\beta$ is specialised to
$0$ then $G_{\theta}(x)$ becomes the Schur polynomial
$s_{\theta}(x)$ so that the values $c_{\theta\mu}^{\nu}(0,0)$
coincide with the Littlewood-Richardson coefficients
$c_{\theta\mu}^{\nu}$ defined by the expansion
$$s_\theta(x) s_\mu(x)=\sum_\nu c_{\theta\mu}^{\nu} s_\nu(x).$$

The coefficients $c_{\lambda\mu}^{\nu}$ with a non-skew partition
$\lambda$ are calculated by the classical Littlewood-Richardson
rule \cite{lr:gc} and its various versions; see e.g. Macdonald
\cite{Mc}, Sagan \cite{sagan}. In the case where $\theta$ is skew,
a rule for calculation of $c_{\theta\mu}^{\nu}$ is given by James
and Peel~\cite{jp:ss} and Zelevinsky~\cite{z:gl} in terms of
combinatorial objects called {\it pictures\/}. There is also a
short proof of a generalised Littlewood-Richardson rule for Schur
polynomials provided by Gasharov \cite{gasharov}, which raises the
question as to whether an analogue exists for Grothendieck
polynomials. A different derivation of such a rule is given by
Molev and Sagan \cite{MS:main}, where a factorial analogue of the
Schur functions was used.

The results given by Buch in \cite{B:buch} are shown to be an
immediate consequence of this new rule. As for the question of
providing a complete description of the Littlewood-Richardson rule
for factorial Grothendieck polynomials, this remains unanswered.

In the last two sections, we turn away from the combinatorial
approach to Grothendieck polynomials used elsewhere in this paper
and consider the so-called double Grothendieck polynomials defined
via isobaric divided difference operators. These chapters work
towards, and eventually prove, the existence of a relationship
between these previously studied double Grothendieck polynomials
and the factorial Grothendieck polynomials introduced here.

\section{Preliminaries}

\subsection{Partitions}

A partition $\lambda=(\lambda_1,\lambda_2,\ldots,\lambda_l)$ is a
finite non-increasing sequence of positive integers,
$\lambda_1\geq\lambda_2\geq\cdots\geq\lambda_l>0$. The number of
parts $l$, is called the length of $\lambda$, and denoted
$\ell(\lambda)$. Throughout this paper, we shall frequently be
dealing with the set of partitions $\lambda$ for which
$\ell(\lambda)\leq n$ for some fixed positive integer $n$. Then,
if $\ell(\lambda)<n$ we shall append zeros to the end of $\lambda$
by defining $\lambda_k=0$ if $\ell(\lambda)<k\leq n$ so we can
treat $\lambda$ as a sequence $\vec{\lambda}$ of $n$ non-negative
integers.

Denote by $|\lambda|$ the weight of the partition $\lambda$,
defined as the sum of its parts,
$|\lambda|=\sum_{i=1}^{\ell(\lambda)}\lambda_i$.

An alternative notation for a partition is to write
$\lambda=(1^{m_1}2^{m_2}\ldots)$ where $m_i$ is the number of
indices $j$ for which $\lambda_j=i$. In such notation, if $m_i=0$
for some $i$, then we omit it from our notation. So for example we
can succinctly write the partition consisting of $n$ parts each
equal to $k$ as $(k^n)$.

The Young diagram of a partition $\lambda$ is formed by
left-aligning $\ell(\lambda)$ rows of boxes, or cells, where the
$i$-th row (counting from the top) contains $\lambda_i$ boxes.



We identify a partition with its Young diagram.

Say $\lambda\to\mu$ if $\mu$ is obtained by adding one box to
$\lambda.$ If $r$ is the row number of the box added to $\lambda$
to create $\mu$ then write $\lambda\stackrel{r}\to\mu$.

By reflecting the diagram of $\lambda$ in the main diagonal, we
get the diagram of another partition, called the conjugate
partition, and denoted $\lambda'$. Alternatively and equivalently,
we can define $\lambda'$ by $\lambda'_j=\#\{i\mid\lambda_i\geq
j\}.$


The main ordering of partitions which we make use of is that of
containment ordering. We say $\lambda\subset\mu$ if the Young
diagram of $\lambda$ is a subset of the Young diagram of $\mu$.

The other ordering which we make mention of is dominance ordering.
We say $\lambda\rhd\mu$ if
$\lambda_1+\cdots+\lambda_k\geq\mu_1+\cdots+\mu_k$ for all $k$.


Suppose we have two partitions $\lambda$, $\mu$ with
$\lambda\supset\mu$. Then we may take the set-theoretic difference
of their Young diagrams and define the skew partition
$\theta=\lambda/\mu$ to be this diagram. Note that every partition
is also a skew partition since $\lambda=\lambda/\phi$ where $\phi$
is the empty partition.

The weight of $\theta$ is the number of boxes it contains:
$|\theta|=|\lambda/\mu|=|\lambda|-|\mu|$.

With regard to notation, the use of $\theta$ shall signify that we
are dealing with a skew partition, while other Greek letters
employed shall refer exclusively to partitions.

\subsection{Tableaux}\label{tableaux}

Let $\theta$ be a skew partition. We introduce a co-ordinate
system of labelling cells of $\theta$ by letting $(i,j)$ be the
intersection of the $i$-th row and the $j$-th column. Define the
content of the cell $\alpha=(i,j)$ to be $c(\alpha)=j-i$.

In each cell $\alpha$ of $\theta$, place a non-empty set
$T(\alpha)\subset[n]=\{1,2,\ldots,n\}$ (a notation we shall use
throughout), such that entries are non-decreasing along rows and
strictly increasing down columns. In other words, if $\alpha$ is
immediately to the left of $\beta$ then $\max(T(\alpha))\leq
\min(T(\beta))$, while if $\alpha$ is immediately above $\beta$,
then $\max(T(\alpha))<\min(T(\beta))$.

An example of such a $(4,4,2,1)/(1)$-tableau is given
in the Introduction.

Such a combinatorial object $T$ is called a semistandard
set-valued $\theta$-tableau. If the meaning is obvious from the
context, we shall often drop the adjectives semistandard and
set-valued. $\th$ is said to be the shape of $T$, which we denote
by ${\rm sh}(T)$.

Define an entry of $T$ to be a pair $(r,\alpha)$ where
$\alpha\in\theta$ is a cell and $r\in T(\alpha)$. Let $|T|$ denote
the number of entries in $T$.

Define an ordering $\prec$ on the entries of $T$ by
$(r,(i,j))\prec(r',(i',j'))$ if $j>j'$, or $j=j'$ and $i<i'$, or
$(i,j)=(i',j')$ and $r>r'$. On occasion, we shall abbreviate this
to $r\prec r'$.

So any two entries of $T$ are comparable under this order, and if
we write all the entries of $T$ in a chain
$(r_1,\alpha_1)\prec(r_2,\alpha_2)\prec\ldots\prec(r_{|T|},\alpha_{|T|})$,
then this is equivalent to reading them one column at a time from
right to left, from top to bottom within each column, and from
largest to smallest in each cell. Writing the entries in this way,
we create a word $r_1r_2\ldots r_{|T|}$, called the {\it column
word\/} of $T$, and denoted $c(T)$.

\subsection{Symmetric functions}

Here we define the monomial symmetric function $m_\lambda$ and the
elementary symmetric function $e_k$ in $n$ variables $\vec{x}$.

For a partition $\lambda=\vec{\lambda}$, define the monomial
symmetric function $m_\lambda$ by
$$m_\lambda(x)=\sum\prod_{i=1}^n x_{\pi(i)}^{\lambda_i}$$
where the sum runs over all distinct values of $\prod_{i=1}^n
x_{\pi(i)}^{\lambda_i}$ that are attainable as $\pi$ runs over the
symmetric group $S_n$.

As an example, if $n=3$, then we have
$m_{(22)}(x_1,x_2,x_3)=x_1^2x_2^2+x_2^2x_3^2+x_3^2x_1^2.$

The elementary symmetric function $e_k$ can now be defined as
$e_k=m_{(1^k)}$.

The monomial symmetric functions $m_\lambda$, where $\lambda$ runs
over all partitions with $\ell(\lambda)\leq n$, form a basis for
the ring of symmetric polynomials in $n$ variables, $\Lambda_n$.

We will stick with convention and use $\Lambda_n$ to denote the
ring of symmetric polynomials in $n$ variables over $\Z$. However,
we will often wish to change the ring of coefficients, so will
often work in $\Lambda_n\otimes_\Z R$ for some ring $R$. As we
shall only ever consider tensor products over $\Z$, the subscript
$\Z$ is to be assumed whenever omitted.

\section{Ordinary Grothendieck Polynomials}

Before starting our work on the factorial Grothendieck
polynomials, first we present some of the theory of the ordinary
Grothendieck polynomials.

\begin{definition} Given a skew diagram $\th$, a field $\F$, $\beta$ an indeterminate over $\F$,
we define the ordinary Grothendieck polynomial
$G_\th(x)\in\F(\beta)[x_1,\ldots ,x_n]$ by
\begin{equation} \label{odef}
G_\th(x)=\sum_T \beta^{|T|-|\theta|} \prod_{\substack{\alpha\in\th
\\ r\in T(\alpha)}} x_r
\end{equation} where the sum is over all semistandard set-valued $\th$-tableau $T$.
\end{definition}

\begin{remark}
In the existing literature, Grothendieck polynomials are often
only presented in the case $\beta=-1$ as a consequence of their
original geometric meaning. The case of arbitrary $\beta$ has been
previously studied in \cite{fkmini} and \cite{FK}, though there is
essentially little difference between the two cases, as can be
seen by replacing $x_i$ with $-x_i/\beta$ in (\ref{odef}) for all
$i$.
\end{remark}


\begin{example} Calculation of $G_{(1)}(x)$.
\end{example}
We can have any nonempty subset of $[n]$ in the single available
cell of $T$, so we have
\[
G_{(1)}(x)=\sum_{\substack{S\subset [n] \\
S\neq\phi}}\beta^{|S|-1}\prod_{i\in S}x_i
=\sum_{j=1}^n\beta^{j-1} \sum_{\substack{S\subset [n] \\
|S|=j}}\prod_{i\in S}x_i =\sum_{j=1}^n \beta^{j-1}e_j(x).
\] where the $e_j$ are the elementary symmetric functions. Hence,
\begin{equation} \label{pi}
1+\beta G_{(1)}(x)=\sum_{j=0}^n\beta^j
e_j(x)=\prod_{i=1}^n(1+\beta x_i)=\Pi(x).
\end{equation}
where for any sequence $y=\vec{y}$, we denote the product
$\prod_{i=1}^n(1+\beta y_i)$ by $\Pi(y)$.

At this stage we will merely state, rather than prove the
following important theorem about ordinary Grothendieck
polynomials, as it is proven in greater generality in Theorems
\ref{symmetric} and \ref{spbasis} of the following section.

\begin{theorem}
The ordinary Grothendieck polynomial $G_{\th}(x)$ is symmetric in
$x_1,\ldots,x_n$, and furthermore the polynomials
$\{G_{\lambda}(x)\mid \ell(\lambda)\leq n\}$ comprise a basis for
the ring of symmetric polynomials in $n$ variables
$\Lambda_n\otimes \F(\beta)$.
\end{theorem}

For a skew-partition $\th$, and partitions $\mu$, $\nu$ with
$\ell(\nu)\leq n$, we define the coefficients
$c_{\theta\mu}^\nu\in \F(\beta)$ by
\begin{equation}\label{deflr}
G_\theta(x)G_\mu(x)=\sum_\nu c^\nu_{\th\mu}G_\nu(x).
\end{equation}
The above theorem shows that these coefficients are well defined.

Before moving onto an important result from the theory of ordinary
Grothendieck polynomials, we present two insertion algorithms
which play an integral role in the proof. Buch \cite{B:buch}
presents a similar column-based insertion algorithm.

First we present a forward row insertion algorithm. As input, this
algorithm takes a set $S\subset[n]$ and a semistandard, set-valued
row $R$ and produces as output a row $R'$ and a set $S'$.

\begin{algorithm}[Forward row insertion algorithm] For all $s\in S$, we perform the following operations
simultaneously:

Place $s$ in the leftmost cell of $R$ such that $s$ is less than
all entries originally in that cell. If such a cell does not
exist, then we add a new cell to the end of $R$ and place $s$ in
this cell.

If there exist entries greater than $s$ occupying cells to the
left of where $s$ was inserted, then remove them from $R$. Call
this a type I ejection. If no such elements exist, then remove
from $R$ all the original entries in the cell $s$ is inserted into
and call this a type II ejection. The resulting row is $R'$ and
the set of elements removed from $R$ is $S'$. \end{algorithm}

For example if $S=\{1,2,3,6,7,8\}$ and $R$ is the row $1, 12, 47,
7, 789, 9$ then the algorithm gives:

$$\begin{array}{cccccccc}
  124678\rightarrow & 1 & 12 & 37 & 7 & 789 & 9 &  \\
  \text{Insert} &  &  & 12 & 46 &  & 78 &  \\
  \text{Eject} &  & 2 & 37 &  & 89 &  &  \\
  \text{Final Result} & 1 & 1 & 12 & 467 & 7 & 789 & \rightarrow23789 \\
\end{array}$$

with output $R'=1, 1, 12, 467, 7, 789$ and $S'=\{2,3,7,8,9\}$.

We show that in this algorithm, if a number $x$ is ejected, then
it is ejected from the rightmost cell in $R$ such that $x$ is
strictly greater than all entries of $R'$ in that cell.

Let $y$ be an entry of $R'$ in the cell $x$ is ejected from, and
suppose that $y\geq x$. If $y$ was, along with $x$ an original
entry of $R$, then $y$ would have been ejected from $R$ at the
same time that $x$ was, a contradiction. Hence $y$ was inserted
from $S$ into $R$. But then, due to the criteria of which cell an
entry gets inserted into, we must have $y<x$, also a
contradiction. So $x$ is greater than all entries of $R'$ in the
cell it was ejected from.

Now consider a cell $\alpha\in R$ to the right of the one $x$ was
ejected from, and let its maximum entry of $\alpha$ in $R'$ be
$y$. If $y$ was an original entry of $R$, then since $R$ is
semistandard, $y\geq x$. Now suppose that $y$ was inserted into
$R$ from $S$, and further suppose, for want of a contradiction,
that $y<x$. Let $z$ be the minimal original entry in $\alpha$. Any
element inserted into $\alpha$ is less than or equal to $y$, so
less than $x$ and hence ejects $x$ via a type I ejection. So no
type II ejections occur in $\alpha$. Now $z>y$ by our insertion
rule for adding $y$, so by maximality of $y$, $z$ must have been
ejected from $R$. Then this must have occurred via a type I
ejection. To be ejected, an element $w<z$ must have been added to
the right of $z$, but such a $w$ cannot be added to the right of
$z$ by the conditions for insertion, a contradiction. Hence $y\geq
x$.

So we have proven that if a number $x$ is ejected, then it is
ejected from the rightmost cell in $R$ such that $x$ is strictly
greater than all entries of $R'$ in that cell. If an element of
$S$, when inserted into $R$ does not cause any entries to be
ejected, then it must have been inserted into a new cell to the
right of $R'$. We are now in a position to describe the inverse to
this algorithm, which we call the reverse row insertion algorithm.

\begin{algorithm}[Reverse Row Insertion Algorithm] The reverse insertion of a set $S'$ into a row $R'$, whose
rightmost cell is possibly denoted special, produces as output a
set $S$ and a row $R$, and is described as follows:

For all $x\in S'$, we perform the following operations
simultaneously:

Insert $x$ in the rightmost non-special cell of $R'$ such that $x$
is strictly greater than all entries already in that cell.

If there exist entries in $R'$ less than $x$ in cells to the
right of $x$, remove them. If this does not occur, then delete all
original entries of $R'$ in the cell in which $x$ was inserted to.

Also, remove all elements in the special cell and delete this
special cell if a special cell exists.

The remaining row is $R$ and $S$ is taken to be the set of all
entries removed from $R'$. \end{algorithm}

We now present an algorithm for inserting a set $S_0\subset[n]$
into a semistandard set-valued tableau $T$.

\begin{algorithm}[Forward Insertion Algorithm]
Let the rows of $T$ be $R_1,R_2,\ldots$ in that order. Step $k$ of
this insertion algorithm consists of inserting $S_{k-1}$ into
$R_k$ using the forward row insertion algorithm described above,
outputting the row $R'_{k}$ and the set $S_k$. The resultant
tableau $T'$ with rows $R'_1,R'_2,\ldots$ is the output of this
algorithm. Write $T'=S\hookrightarrow T$.
\end{algorithm}

Now we show that $T'=S\hookrightarrow T$ is a semistandard
set-valued tableau, and furthermore that if $T$ has shape
$\lambda$ and $T'$ has shape $\mu$, then $\lambda\tangle\mu$.

It is an immediate consequence of the nature of the row insertion
algorithm that each row of $T'$ is non-decreasing. To show that
entries strictly decrease down a column, we need to look at what
happens to an entry ejected from a row $R_k$ and inserted into
$R_{k+1}$. Suppose that this entry is $a$ and is ejected from the
$j$'th column and inserted into the $i$'th column of $R_{k+1}$.
Then $a\in T(k,j)$ so $a<T(k+1,j)$ and hence $i\leq j$. Any entry
in $T(k,i)$ greater than or equal to $a$ must also be ejected from
$R_k$ so $T'(k,i)<a$. Since this algorithm always decreases the
entries in any given cell, the only place where semistandardness
down a column needs to be checked is of the form $T'(k,i)$ above
the inserted $a$ as checked above, so $T'$ is indeed semistandard.

In the transition from $T$ to $T'$, clearly no two boxes can be
added in the same row. Now, when a box is added, no entries are
ejected from this box. We have just shown above that the path of
inserted and ejected entries always moves downward and to the
left, so it is impossible for entries to be added strictly below
an added box, so hence no two boxes can be added in the same
column, so our desired statement regarding the relative shapes of
$T$ and $T'$ is proven.

We now construct the inverse algorithm. Let $\lambda$ be a
partition and suppose $T'$ is a semistandard set-valued tableau
with shape $\mu$ where $\lambda\tangle\mu$. Call a cell of $T'$
special if it is in $\mu/\lambda$. The inverse algorithm takes as
input $T'$ as described above and produces a $\lambda$-tableau $T$
and a set $S\subset[n]$ for which $T'=S\hookrightarrow T$.

 Now,
supposing we have a $\mu$-tableau $T'$ as described above with
rows $R'_1,R'_2,\ldots,R'_{\ell(\mu)}$. Let $S_{\ell(\mu)}=\phi$
and form $R_k$ and $S_{k-1}$ by reverse inserting $S_k$ into
$R'_k$. Then $T$ is the resulting tableau consisting of rows
$R_1,R_2,\ldots$ and $S=S_0$.

This completes our description of the necessary insertion
algorithms. We note that the forward row insertion algorithm and
the reverse row insertion algorithm are inverses of each other, we
have constructed the inverse of the map $(S,T)\mapsto
(S\hookrightarrow T)$ and hence this map is a bijection.

The following equation is due to Lenart \cite{Lenart}. The proof
we give however is based on the algorithm depicted above.

Say $\lambda\tangle\mu$ if $\mu/\lambda$ has all its boxes in
different rows and columns (this notation also includes the case
$\lambda=\mu$). If we want to discount the possibility that
$\lambda=\mu$, then we write $\lambda\tangle^*\mu$.

\begin{proposition} \label{buch} \cite{Lenart}
\[
G_\lambda(x)\Pi(x)=\sum_{\lambda\tangle\mu}\beta^{|\mu/\lambda|}G_\mu(x).
\]
\end{proposition}

\begin{proof} We have a bijection via our insertion algorithm between pairs $(S,T)$ with $S\subset [n]$
and $T$ a $\lambda$-tableau, and $\mu$-tableau $T'$ where $\mu$ is
a partition such that $\lambda\tangle\mu$. Furthermore, if we let
$x^T=\prod_{r\in T}x_r$, we note that the insertion algorithm at
no time creates destroys or changes the numbers occurring in the
tableau, only moves them and thus $x^T x^S=x^{(S\hookrightarrow
T)}$.

Therefore,
  \begin{eqnarray*}
G_\lambda(x)\Pi(x)&=&\sum_{\rm{sh}(T)=\lambda}\beta^{|T|-|\lambda|}x^T
\sum_{S\subset[n]}\beta^{|S|}x^S \\
&=&\sum_{(T,S)}\beta^{|T|+|S|-|\lambda|}x^{(S\hookrightarrow T)}
\\
\hspace{30pt}&=&\sum_{\lambda\tangle\mu}\beta^{|\mu/\lambda|}\sum_{\rm{sh}(T')=\mu}
\beta^{|T'|-|\mu|}x^{T'} \\
&=&\sum_{\lambda\tangle\mu}\beta^{|\mu/\lambda|}G_\mu(x)
\end{eqnarray*}
as required. \end{proof}

This last result provides the values of $c_{\lambda(1)}^\nu$ for
all partitions $\lambda$ and $\nu$. Later, we shall prove Theorem
\ref{main} providing a rule describing the general coefficient
$c_{\th\mu}^{\nu}$. This theorem encompasses two special cases
which are known thanks to Buch \cite{B:buch}, namely that when
$\th$ is a partition, and when $\mu=\phi$, the empty partition. We
shall finish off this section by quoting these results. In order
to do so however, we first need to introduce the idea of a lattice
word.

\begin{definition}
We say that a sequence of positive integers
$w=(i_1,i_2,\ldots,i_l)$ has content $(c_1,c_2,\ldots)$ if $c_j$
is equal to the number of occurrences of $j$ in $w$. We call $w$ a {\it lattice word} if for each $k$, the content of the
subsequence $(i_1,i_2,\ldots ,i_k)$ is a partition.
\end{definition}

For the case where $\theta=\lambda$, a partition, Buch's result is
as follows:
\begin{theorem} \label{b1} \cite{B:buch}
$c_{\lambda\mu}^{\nu}$ is equal to $\beta^{|\nu|-|\lambda|-|\mu|}$
times the number of set-valued tableaux $T$ of shape $\lambda*\mu$
such that $c(T)$ is a lattice word with content $\nu$.
\end{theorem}
Here, $\lambda*\mu$ is defined to be the skew diagram obtained by
adjoining the top right hand corner of $\lambda$ to the bottom
left corner of $\mu$ as shown in the diagram below.

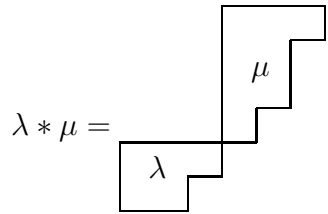
\setlength{\unitlength}{1.1em}
$$\lambda*\mu=
\parbox[c]{6\unitlength}{\begin{picture}(6,6)
\put(3,6){\line(1,0){3}} \put(5,5){\line(1,0){1}}
\put(4,3){\line(1,0){1}} \put(0,2){\line(1,0){4}}
\put(2,1){\line(1,0){1}} \put(0,0){\line(1,0){2}}

\put(0,0){\line(0,1){2}} \put(2,1){\line(0,-1){1}}
\put(3,6){\line(0,-1){5}}\put(4,3){\line(0,-1){1}}\put(5,5){\line(0,-1){2}}
\put(6,6){\line(0,-1){1}}

\put(0.84,0.95){$\lambda$} \put(3.84,3.96){$\mu$}

\end{picture}
\setlength{\unitlength}{1pt} }
$$

For the case where $\mu=\phi$, the empty partition, Buch's result,
expanding the skew Grothendieck polynomial $G_\theta(x)$ in the
basis $\{ G_{\lambda}(x)\mid \ell(\lambda)\leq n \}$ is as
follows:
\begin{theorem} \label{b2}\cite{B:buch}
$c_{\th\phi}^\nu$ is equal to the number of set-valued tableaux of
shape $\th$ such that $c(T)$ is a lattice word with content $\nu$.
\end{theorem}

\section{The Factorial Grothendieck Polynomials}

Now we are ready to begin our study of the factorial Grothendieck
polynomials, the main focus of this paper. Again, we work over an
arbitrary field $\F$, and let $\beta$ be an indeterminate over $\F$.
In addition to this, we shall also have to introduce a second family
of variables as part of the factorial Grothendieck polynomials.

Define the binary operation $\oplus$ (borrowed from \cite{fkmini}
and \cite{FK}) by
$$x\oplus y=x+y+\beta xy$$
and denote the inverse of $\oplus$ by $\ominus$, so we have
$\ominus x=\frac{-x}{1+\beta x}$ and $x\ominus
y=\frac{x-y}{1+\beta y}$.

\subsection{Definition and basic properties}

Let $\theta$ be a skew diagram, $a=(a_k)_{k\in\Z}$ be a sequence
of variables (in the most important case, where $\theta$ is a
partition, we only need to consider $(a_k)_{k=1}^{\infty}$). We
are now in a position to define the factorial Grothendieck
polynomials in $n$ variables $x=(x_1,x_2,\ldots,x_n)$.

\begin{definition}[Factorial Grothendieck Polynomials] The factorial
Grothendieck polynomial $G_\theta(x|a)$ is defined to be
\begin{equation} \label{eq:def}
G_\theta(x|a)=\sum_T \beta^{|T|-|\theta|}
\prod_{\substack{\alpha\in
\theta \\
r\in T(\alpha)}}x_r\oplus a_{r+c(\alpha)},
\end{equation}
recalling that $c(\alpha)$ is the content of the cell $\alpha$,
defined by $c(i,j)=j-i$. The summation is taken over all
semistandard set-valued $\theta$-tableaux $T$.
\end{definition}

\underline{Remarks}

1. The name {\it factorial Grothendieck polynomial} is chosen to
stress the analogy with the factorial Schur functions, as
mentioned for example (though not explicitly with this name), in
variation 6 of MacDonald's theme and variations of Schur functions
\cite{Mc:Schur}. The factorial Schur functions are obtainable as a
specialisation of the factorial Grothendieck polynomials by
setting $\beta=0$, though to be truly consistent with the
established literature, one should accompany this specialisation
with the transformation $a\mapsto -a$.

2. Setting $\theta=\phi$, the empty partition, we get
$G_\phi(x|a)=1$.

3. The $\beta$ can be seen to play the role of marking the degree.
For if we assign a degree of $-1$ to $\beta$, where $x$ and $a$
each have degree 1, then $G_\theta(x|a)$ becomes homogenous of
degree $|\theta|$.

4. If we set $a=0$, then we recover the ordinary Grothendieck
polynomials through specialisation.

5. If $\lambda$ is a partition with $\ell(\lambda)>n$, then it is
impossible to fill the first column of $\lambda$ to form a
semistandard tableau, so $G_\lambda(x|a)=0$. Hence we tend to work
only with partitions of at most $n$ parts.

6. In a similar vein to the connection between factorial Schur
functions and double Schubert polynomials, as pointed out by Lascoux
\cite{ldouble}, there exists a relationship between these factorial
Grothendieck polynomials and the double Grothendieck polynomials
discussed for example in \cite{B:buch}, amongst other places. The
final two sections of this paper work towards proving such a result,
culminating in Theorem \ref{final}, which provides a succinct
relationship between these two different types of Grothendieck
polynomials.

\begin{example} Let us calculate $G_{(1)}(x|a)$. Here we use
$x\oplus a$ to represent the sequence $(x_1\oplus a_1,x_2\oplus
a_2,\ldots,x_n\oplus a_n)$.

Similarly to the calculation of $G_{(1)}(x)$, we can have any
nonempty subset of $[n]$ in the one box of $T$, so we have

\begin{eqnarray*}
G_{(1)}(x|a)&=&\sum_{\substack{S\subset [n] \\
S\neq\phi}}\beta^{|S|-1}\prod_{i\in S}x_i\oplus a_i
=\sum_{j=1}^n\beta^{j-1} \sum_{\substack{S\subset [n] \\
|S|=j}}\prod_{i\in S}x_i\oplus a_i \\
&=&\sum_{j=1}^n \beta^{j-1}e_j(x\oplus a),
\end{eqnarray*} where the $e_j$ are the elementary symmetric functions. Hence,\[
1+\beta G_{(1)}(x|a)=\sum_{j=0}^n\beta^j e_j(x\oplus
a)=\prod_{j=1}^n (1+\beta (x_j\oplus a_j))=\Pi(x)\Pi(a).
\]

\end{example}

\begin{theorem} \label{symmetric} The factorial Grothendieck polynomials are
symmetric in $x_1,x_2,\ldots,x_n$.
\end{theorem}

\begin{proof} (This proof is a generalisation of a standard argument,
for example as appears in \cite[Prop 4.4.2]{sagan}.) The symmetric
group $S_n$ acts on the ring of polynomials in $n$ variables
$x_1,x_2,\ldots,x_n$ by permuting variables: $\pi
P(x_1,\ldots,x_n)=P(x_{\pi(1)},\ldots,x_{\pi(n)})$ for $\pi\in
S_n$. Since the adjacent transpositions $(i,i+1)$ generate $S_n$,
to show that $G_\theta(x|a)$ is symmetric it suffices to show that
$G_\theta(x|a)$ is stable under interchanging $x_i$ and $x_{i+1}$.
We consider marked semistandard tableaux, with an entry $j$ marked
in one of three ways:

\begin{enumerate}
\item first marking - $j$ - corresponding to taking the $x$ term
from $x\oplus a$. \item second marking - $j^*$ - taking the $a$
term from $x\oplus a$. \item third marking - $jj^*$ - taking the
$\beta xa$ term from $x\oplus a$.
\end{enumerate}

Then we can write $G_\theta(x|a)$ as a sum over marked tableaux,
where each marked tableau $T$ contributes the monomial

$$w(T)=\beta^{|T|-|\theta|}\prod_{r \text{\ unstarred}}x_r\prod_{r
\text{\ starred}}a_{r+c(\alpha)}$$ to the sum
$G_\theta(x|a)=\sum_T w(T)$.


\[
\setlength{\unitlength}{1.6em}
\begin{picture}(3,3)
\multiput(0,0)(0,1){3}{\line (1,0){2}}
\multiput(0,0)(1,0){3}{\line(0,1){2}}
\multiput(1,3)(0,-1){2}{\line (1,0){2}} \multiput(1,2
)(1,0){3}{\line(0,1){1}} \put(0.1,0.3){$2^*6$}
\put(0.3,1.3){$1^*$} \put(1.2,0.3){$78$}
\put(1.0,1.3){$455^{\!*}$} \put(1.2,2.3){$13$}
\put(2.1,2.3){$4^*9$}
\end{picture}
\setlength{\unitlength}{1pt}
\]
As an example, if $T$ is the above tableau, then we have
$$w(T)=\beta^7x_9a_6x_3a_2x_5a_5x_4x_8x_7a_0x_6a_0.$$

Note that there is no ambiguity between first and third, or second
and third markings, since the same number cannot occur twice in
the same cell.

We now find a bijection $T\rightarrow T'$ between marked tableaux
such that $(i,i+1)w(T)=w(T')$.

Given $T$, we construct $T'$ as follows:

All entries not $i$ or $i+1$ remain unchanged.

If $i$ and $i+1$ appear in the same column, we swap their
markings. An example with $i=2$ is the following:


\setlength{\unitlength}{1.6em}
$$\begin{picture}(1,2)
\multiput(0,0)(1,0){2}{\line(0,1){2}}
\multiput(0,0)(0,1){3}{\line(1,0){1}} \put(0.18,1.3){$22^{*}$}
\put(0.36,0.3){$3$}
\end{picture}
\quad\raisebox{1.34em}{$\longleftrightarrow$}\quad
\begin{picture}(1,2)
\multiput(0,0)(1,0){2}{\line(0,1){2}}
\multiput(0,0)(0,1){3}{\line(1,0){1}} \put(0.36,1.3){$2$}
\put(0.18,0.3){$33^*$}
\end{picture}
\setlength{\unitlength}{1pt}
$$

All other occurrences of $i$ and $i+1$ are called free, and we
deal with them one row at a time, independently of each other.

Suppose that there are $a$ free $i$'s and $b$ free $i+1$'s in a
row. Here we are not counting starred markings and also not
distinguishing between an unstarred number in the first marking
and in the third marking. If $a=b$, the row remains unchanged. Now
let us assume that $a>b$. (The $a<b$ case can proceed similarly,
or alternatively and equivalently can be defined to be the inverse
of the $a>b$ case.)

Consider those cells from the $(b+1)$-th free $i$ to the $a$-th
free $i$ inclusive, and call these cells $L$. If the rightmost box
of $L$ contains an unstarred $i+1$, we extend $L$ to the left to
start at the $b$-th free $i$. We leave boxes outside of $L$
unchanged and modify boxes in $L\subset T$ to form $L'\subset T'$
as follows:

1. For each second marking $i^*$ in $L$, not in the leftmost box
of $L$, replace it by an $(i+1)^*$ in $L'$ one box to the left.
Similarly, for each third marking $ii^*$ in $L$, not in the
leftmost box of $L$, we replace it by an $(i+1)(i+1)^*$ in $L'$
one box to the left.

2. Place an $i+1$ in any cells of $L'$ which are still empty.

3. Any $i^*$ in the leftmost box of $L$ or an $(i+1)^*$ in the
rightmost box of $L'$ is left unchanged.

4. If there is an $i+1$ in the last square of $L$, place an $i$ in
the first square of $L'$.

To illustrate this more clearly, we provide now an example of the
bijection between free rows (note in this example, $i=2$, $a=3$,
$b=1$, and $L$ runs from the third to the sixth cells in the row
inclusive).

\setlength{\unitlength}{1.6em}
$$\parbox[c]{7\unitlength}{\begin{picture}(7,1)
\multiput(0,0)(1,0){8}{\line(0,1){1}}
\multiput(0,0)(0,1){2}{\line(1,0){7}} \put(0.35,0.3){$2^*$}
\put(1.36,0.3){$2$} \put(2.18,0.3){$22^{*}$} \put(3.35,0.3){$2^*$}
\put(4.36,0.3){$2$} \put(5.18,0.3){$23^*$} \put(6.36,0.3){$3$}
\end{picture}}
\quad\longleftrightarrow\quad
\parbox[c]{7\unitlength}{\begin{picture}(7,1)
\multiput(0,0)(1,0){8}{\line(0,1){1}}
\multiput(0,0)(0,1){2}{\line(1,0){7}} \put(0.35,0.3){$2^*$}
\put(1.36,0.3){$2$} \put(2.01,0.3){$2^{*}\hspace{-1pt}3^{*}$}
\put(3.36,0.3){$3$} \put(4.36,0.3){$3$} \put(5.18,0.3){$33^*$}
\put(6.36,0.3){$3$}
\end{picture}}
\setlength{\unitlength}{1pt}
$$

From the structure of the construction of the map $T\mapsto T'$,
we can easily see that it is an involution, so thus is bijective,
and furthermore except that the number of $i$'s and the number of
$(i+1)$'s is reversed, weights are preserved, in the sense that we
have the desired equation $(i,i+1)w(T)=w(T')$. Thus we have

\[
(i,i+1)G_\theta(x|a)=\sum_T (i,i+1)w(T)=\sum_{T'}
w(T')=G_\theta(x|a)
\]
as required, so the proof is complete.
\end{proof}

Given a partition
$\lambda=(\lambda_1,\lambda_1,\ldots,\lambda_n)$, define the
sequence $a_\lambda=\vec{(a_\lambda)}$ by $(a_\lambda)_i= \ominus a_{n+1-i+\lambda_i} =\frac{-a_{n+1-i+\lambda_i}}{1+\beta
a_{n+1-i+\lambda_i}}$

\begin{theorem}[Vanishing Theorem]

Suppose $\lambda$ and $\mu$ are partitions with $\ell(\lambda)\leq
n$. Then
$$G_{\lambda}(a_\mu|a)=0\quad\hbox{if}\quad\lambda\not\subset\mu$$
$$G_{\lambda}(a_\lambda|a)\neq0$$
\end{theorem}

\begin{proof} (This argument is borrowed from Okounkov's
paper~\cite{M:vanish} and is included here for completeness and
its importance)

Since $G_\lambda(x|a)$ is symmetric in $x$, we replace the
sequence $\vec{x}$ by $(x_n,x_{n-1},\ldots,x_1)$ in (\ref{eq:def})
to obtain

\[
G_\lambda(x|a)=\sum_T\beta^{|T|-|\lambda|}\prod_{\substack{\alpha\in\lambda \\
r\in T(\alpha)}}x_{n+1-r}\oplus a_{r+c(\alpha)}.
\]

Thus we have

\begin{eqnarray}\label{xtot}
G_\lambda(a_\mu|a)&=&\sum_T\beta^{|T|-|\lambda|}x^T \\
\text{where}\qquad x^T&=&\prod_{\substack{\alpha\in\lambda \\
r\in T(\alpha)}}\frac{a_{r+c(\alpha)}-a_{r+\mu_{n+1-r}}}{1+\beta
a_{r+\mu_{n+1-r}}}.\nonumber
\end{eqnarray}

In order to continue, we need the following proposition.

\begin{proposition}\label{prop}
$x^T\neq 0$ if and only if $T(i,j)\geq n+i-\mu'_j$ for all
$(i,j)\in \lambda$.
\end{proposition}

\begin{proof} (As well as representing a set, sometimes we write
$T(\alpha)$ here and treat it like an integer, to do so means that
the result holds for any element of $T(\alpha)$.)

We have
\begin{eqnarray} \label{vprop}
x^T&\neq& 0\nonumber\\
\iff \hspace{0pt}\mu_{n+1-T(\alpha)}&\neq&
c(\alpha)\quad\forall\alpha
\end{eqnarray}

Assuming that this holds, we shall now prove by induction on $j$
that $\mu_{n+1-T(1,j)}\geq j$.

For $j=1$, by (\ref{vprop}), $\mu_{n+1-T(1,1)}\neq 0$ so
$\mu_{n+1-T(1,1)}\geq 1$ as it is a non-negative integer.

Now suppose $j>1$ and $\mu_{n+1-T(1,j-1)}\geq j-1$.

As $T(1,j)\geq T(1,j-1)$,
$\mu_{n+1-T(1,j)}\geq\mu_{n+1-T(1,j-1)}$.

But by (\ref{vprop}), $\mu_{n+1-T(1,j)}\neq j-1$ so it must be
that $\mu_{n+1-T(1,j)}\geq j$ as required and the induction is
complete.

Therefore $n+1-T(1,j)\leq\mu_j'$, i.e. $T(1,j)\geq n+1-\mu'_j$.
Since $T(i+1,j)>T(i,j)$, a straightforward induction on $i$ gives
$T(i,j)\geq n+i-\mu_j'$.

Now suppose that $T(i,j)\geq n+i-\mu'_j$ for all $(i,j)\in
\lambda$. Equivalently this can be written as $n+i-T(i,j)\leq
\mu_j'$, which gives us the chain of inequalities
$\mu_{n+1-T(i,j)}\geq\mu_{n+i-T(i,j)}\geq\mu_{\mu'_j}\geq j>j-i$.
In particular, this shows $\mu_{n+1-T(\alpha)}\neq c(\alpha)$, so
$x^T\neq 0$ as required and the proof of the proposition is
complete.
\end{proof}

Now we return to proving the vanishing theorem and apply the
condition that $T(i,j)\leq n$ to the above proposition. If
$G_\lambda(a_\mu|a)\neq 0$, then we have:

\[
\lambda_j'\geq i\implies (i,j)\in \lambda\implies n+i-\mu'_j\leq
n\implies \mu_j'\geq i.
\]

Thus $\lambda_j'\leq\mu_j'$ for all $j$, so $\lambda\subset\mu$
and the first part of the vanishing theorem is proven.

In the case $\lambda=\mu$, equality must hold everywhere, so there
is a unique $\lambda$-tableau $T$ for which $x^T\neq 0$, namely
that with $T(i,j)=n+i-\lambda_j'$ for all $(i,j)\in T$. Hence
$G_\lambda(a_\lambda|a)\neq 0$.
\end{proof}

From the above, we can write an explicit formula for
$G_\lambda(a_\lambda|a)$. After making the change $i\to
1+\lambda_j'-i$ to neaten up the result, we get:

$$G_\lambda(a_\lambda|a)=\prod_{(i,j)\in\lambda}\frac{a_{n+j-\lambda'_j}-a_{\lambda_i+n-i+1}}{1+\beta a_{\lambda_i+n-i+1}}.$$

We pause to introduce a space utilised in a couple of subsequent
proofs. Let $L_k$ denote the subspace of $\Lambda_n\otimes
\F(\beta,a)$ spanned by the monomial symmetric functions
$\{m_\lambda\mid\lambda\subset k^n\}$.

\begin{theorem} \label{basis} $\{G_\lambda(x|a)\}$ form a basis in $\Lambda_n\otimes_\Z \F(\beta,a)$ as
$\lambda$ runs over all partitions with $\ell(\lambda)\leq n$.
\end{theorem}

\begin{proof} If $\mu\subset k^n$, then $G_\mu(x|a)\in L_k$, for a
number $i$ in a $\mu$-tableau $T$ can appear at most once in each
column, and hence at most $k$ times overall so the exponent of
$x_i$ in $G_\mu(x|a)$ is at most $k$.

Let $\rho_1<\rho_2<\cdots<\rho_{l}$ be a fixed ordering of the
$l={{n+k}\choose{k}}$ partitions contained in $(k^n)$ which is a
refinement of the dominance ordering. Define the matrix $D_k$ by
\begin{equation} \label{matrix}
\left(%
\begin{array}{c}
  G_{\rho_1}(x|a) \\
  \vdots \\
  G_{\rho_l}(x|a) \\
\end{array}%
\right)=D_k\left(%
\begin{array}{c}
  m_{\rho_1}(x) \\

  \vdots \\
  m_{\rho_l}(x) \\
\end{array}%
\right) \end{equation} which is possible since
$\{m_\lambda\mid\lambda\subset k^n\}$ forms a basis of $L_k$.

Let $d_k=d_k(\beta,a)=\det{D_k}$.

From the definition of $G_\lambda(x|a)$, we see that every entry
of $D_k$ is a polynomial in $\beta$ and $a$ and hence the same is
true of $d_k$.

If we specialise to the case $\beta=0$, $a=0$, then $D_k$ becomes
the transition matrix from the monomial symmetric functions to the
classical Schur functions. In \cite[Ch 1, 6.5]{Mc}, this
transition matrix is shown to be lower triangular, with 1's along
the main diagonal, so has determinant 1 and thus $d_k(0,0)=1$.
Hence $d(\beta,a)$ is not identically zero, so $D_k$ is invertible
and thus $\{G_\lambda(x|a)\mid \lambda\subset k^n\}$ is a basis of
$L_k$.

As $L_0\subset L_1\subset L_2\subset\cdots$ and $\cup_{k=0}^\infty
L_k=\Lambda_n\otimes \F(\beta,a)$,
$\{G_\lambda(x|a)\mid\ell(\lambda)\leq n\}$ forms a basis for
$\Lambda_n\otimes\F(\beta,a)$. \end{proof}

\subsection{Analysis of poles}

In this section, we do not make any use of skew diagrams, so only
need to deal with the sequence of variables
$(a_k)_{k=1}^{\infty}$.

Suppose $P(x)\in\Lambda_n\otimes\F[\beta,a]$, and suppose we
expand $P(x)$ in the basis of factorial Grothendieck polynomials
$G_\lambda(x|a)$.

\begin{equation} \label{expand}
P(x)=\sum_\lambda d_\lambda G_\lambda(x|a)
\end{equation}

The coefficients $d_\lambda$ can be written as a quotient of
coprime polynomials in $\beta$ and $a$:
$d_\lambda=f_\lambda/g_\lambda$.

\begin{lemma} \label{polecontrol}
 The only possible irreducible factors of $g_\lambda$ are of the form $1+\beta a_i$ for some $i>0$.
\end{lemma}

\begin{proof} First, fix a $k$ such that $P(x)\in L_k$.
If we first expand $P(x)$ in the basis of monomial symmetric
functions $m_\lambda$, we find that the coefficients are all
polynomials in $\beta$ and $a$. Using (\ref{matrix}) and Cramer's
rule to subsequently determine the $d_\lambda$, we find that the
denominators $g_\lambda$ must all divide $d_k$.

Setting $x=a_\mu$ in (\ref{expand}) and applying the vanishing
theorem gives
$$d_\mu=\frac{1}{G_\mu(a_\mu|a)}\biggl(P(a_\mu)-\sum_{\rho\subsetneq\mu}d_\rho
G_\rho(a_\mu|a)\biggr).$$

This provides a recurrence from which the coefficients $d_\lambda$
can be computed inductively using inclusion ordering. From such an
induction, we can conclude that the only possible irreducible
factors of the denominators $g_\lambda$ are of the form $1+\beta
a_i$ (from poles of $G_\rho(a_\mu|a)$) or $a_i-a_j$ (from zeros of
$G_\mu(a_\mu|a)$). If the latter occurs, then $g_\lambda(0,0)=0$,
contradicting $d_k(0,0)=1$ as $g_\lambda|d_k$. Thus the only
possible irreducible factors of denominators $g_\lambda$ are of
the form $1+\beta a_i$ ($i>0$). \end{proof}

In fact we can prove that the only possible irreducible factors of
$d_k(\beta,a)$ are of the form $1+\beta a_i$ for some $i$. For if
$f$ is irreducible and $f$ divides $d_k$, then working over the
integral domain $\F[\beta][a]/(f)$, where $(f)$ is the ideal
generated by $f$, we have that the determinant of the transition
matrix from the monomial symmetric functions to the factorial
Grothendieck polynomials is zero. Hence the factorial Grothendieck
polynomials are linearly dependent. So there exist
$c_\lambda\in\F[\beta][a]/(f)$ not all zero such that
$\sum_\lambda c_\lambda G_\lambda(x|a)=0.$ If
$b_\lambda\in\F[\beta][a]$ is such that $c_\lambda=b_\lambda+(f)$
then $\sum_\lambda b_\lambda G_\lambda(x|a)=fg$ for some
$g\in\F[\beta][a][x]$ and not all $b_\lambda$ are divisible by
$f$. Then $g=\sum_\lambda \frac{b_\lambda}{f} G_\lambda$ and since
not all $b_\lambda$ are divisible by $f$, $f$ can appear as a
denominator of an expansion of the form (\ref{expand}), and hence
from our above result concerning such denominators, $f$ must be of
the form $1+\beta a_i$ for some $i$.

In the subsequent section, we shall prove the following formula,
which shows that for all $i>0$, $1+\beta a_i$ can appear as a
factor of a denominator in an expansion of the form
(\ref{expand}), and hence divides $d_k$ for large enough $k$.

\begin{proposition} \label{main1prop}
\begin{equation} \label{main1}
G_\lambda(x|a)\Pi(x)=\Pi(a_\lambda)\sum_{\lambda\tangle\mu}\beta^{|\mu/\lambda|}G_\mu(x|a)
\end{equation}
\end{proposition}

Once this formula is proven, we have the stronger result.

\begin{theorem} \label{spbasis}
The specialisation of $\{G_\lambda(x|a)\mid\ell(\lambda)\leq n\}$ under an evaluation
homomorphism $\F[\beta,a]\rightarrow\F$ forms a basis of
$\Lambda_n\otimes\F$ if and only if
$a_i\beta\neq-1$ for all $i$.
\end{theorem}

Note that this also includes the important case of the ordinary
Grothendieck polynomials via the specialisation $a=0$.

\section{A Recurrence for the Coefficients}

\subsection{Proof of Proposition \ref{main1prop}}

Define coefficients $c^\mu_\lambda=c^\mu_\lambda(\beta,a)$ by
\begin{equation}\label{cmndef}\frac{G_\lambda(x|a)\Pi(x)}{\Pi(a_\lambda)}=\sum_{\mu}\beta^{|\mu|-|\lambda|}c^\mu_\lambda
G_\mu(x|a).\end{equation} These are well defined since the
factorial Grothendieck polynomials are known to form a basis.
(Theorem \ref{basis}.)

First, consider (\ref{cmndef}) with $x$ and $a$ replaced by
$x/\beta$ and $a/\beta$ respectively:
\[
\frac{\beta^{|\lambda|}G_\lambda(\frac{x}{\beta}|\frac{a}{\beta})\Pi(\frac{x}{\beta})}
{\Pi(\frac{a_\lambda}{\beta})}=\sum_{\mu}c^\mu_\lambda\Big(\beta,\frac{a}{\beta}\Big)\beta^{|\mu|}
G_\mu\Big(\frac{x}{\beta}|\frac{a}{\beta}\Big).
\]
Terms of the form $\beta^{|\nu|}
G_\nu(\frac{x}{\beta}|\frac{a}{\beta})$ and $\Pi(\frac{y}{\beta})$
are both independent of $\beta$. Hence
$c_\lambda^\mu(\beta,\frac{a}{\beta})$ is also independent of
$\beta$. As we already know $c_\lambda^\mu$ is a rational function
of $\beta$ and $a$, this last piece of information tells us that
in fact $c_\lambda^\mu$ is a rational function of $\beta a_1,\beta
a_2,\ldots$.

Setting $x=a_\mu$ in (\ref{cmndef}) and applying the vanishing
theorem gives:

\begin{equation} \label{pole}
c^\mu_\lambda=\frac{1}{\beta^{|\mu|-|\lambda|}G_\mu(a_\mu|a)}
\biggl(\frac{G_\lambda(a_\mu|a)\Pi(a_\mu)}{\Pi(a_\lambda)}-\sum_{\rho\subsetneq\mu}\beta^{|\rho|-|\lambda|}c^\rho_\lambda
G_\rho(a_\mu|a)\biggr)
\end{equation}

from which we compute the coefficients $c_\lambda^\mu$ inductively
on $\mu$.

If $\rho$ is a minimal partition with respect to containment order
for which $c^\rho_\lambda\neq 0$, then this gives
$G_\lambda(a_\rho|a)\neq 0$, so by the vanishing theorem,
$\lambda\subset\rho$. So we may rewrite our sum in (\ref{pole}) as
a sum over $\lambda\subset\rho\subsetneq\mu$.

We shall now prove by induction on $\mu$ that
$c_\lambda^\mu\in\F[\beta,a]$. So suppose that
$c_\lambda^\rho\in\F[\beta,a]$ for all $\rho\subsetneq\mu$.

From (\ref{pole}), we find that the following list gives all
possibilities for poles of $c_\lambda^\mu$:
\begin{enumerate}
\item zeros of $\beta^{|\mu|-|\lambda|}G_\mu(a_\mu|a)$. \item
poles of $\beta^{|\rho|-|\lambda|} c_\lambda^\rho
G_\rho(a_\mu|a)$, where $\lambda\subset\rho \subsetneq \mu$. \item
poles of $G_\lambda(a_\mu|a)\Pi(a_\mu)\Pi(a_\lambda)^{-1}$.
\end{enumerate}

1. Zeros of $\beta^{|\mu|-|\lambda|}G_\mu(a_\mu|a)$ are of the
form $\beta$ or $a_i-a_j$. However
$G_\lambda(x|a)\Pi(x)\Pi(a_\lambda)^{-1}\in\Lambda_n\otimes\F[\beta,a]$ since $\Pi(a_\lambda)^{-1}=\prod_{i=1}^n (1+\beta a_{n+1-i+\lambda_i})$,
so by Lemma \ref{polecontrol}, the poles of
$\beta^{|\mu|-|\lambda|}c_\lambda^\mu$ can only have irreducible
factors of the form $1+\beta a_i$. This leaves open the
possibility that $\beta$ could be a pole of $c_\lambda^\mu$. If
this were the case, since $c_\lambda^\mu$ is a rational function
of $\beta a$, $c_\lambda^\mu$ would also have to contain a pole
which vanishes at $a=0$, contradicting our general result
concerning poles.

2. A pole of $\beta^{|\rho|-|\lambda|}$ cannot be a pole of
$c_\lambda^\mu$ by the argument above. By inductive assumption,
there do not exist any poles of $c_\lambda^\rho$, where $\rho
\subsetneq\mu$. Now write $G_\rho(a_\mu|a)=\sum_T
\beta^{|T|-|\rho|}x^T$ as per (\ref{xtot}). By our proposition, if
$x^T\neq 0$, then $T$ can have at most the entries
$n+1-\mu_j',n+2-\mu_j',\ldots,n$ in the $j$-th column. This
maximal set of entries are exactly the entries of the unique
tableau $T$ which contributes a non-zero amount to
$G_\mu(a_\mu|a)$. Hence the pole of $x^T$ is at most that of
$G_\mu(a_\mu|a)$, so since we divide by $G_\mu(a_\mu|a)$, this
gives no contributions to poles of $c_\lambda^\mu$.

3. Write $G_\lambda(a_\mu|a)=\sum_T \beta^{|T|-|\lambda|}x^T$ as
per (\ref{xtot}). Then we have
$$G_\lambda(a_\mu|a)\frac{\Pi(a_\mu)}{\Pi(a_\lambda)}=\sum_T
\beta^{|T|-|\lambda|}x^{T}\frac{\Pi(a_\mu)}{\Pi(a_\lambda)}.$$
Suppose that $c_\lambda^\mu$ has a factor $1+\beta a_k$ in its
denominator. Then $1+\beta a_k$ is a pole of $x^T$ or $\Pi(a_\mu)$
so is of the form $1+\beta a_{i+\mu_{n+1-i}}$ for some $i\in[n]$.

We only need to consider those tableaux $T$ for which the
multiplicity of the factor $1+\beta a_{i+\mu_{n+1-i}}$ in the
denominator of $x^T\Pi(a_\mu)$ is strictly greater than the
corresponding multiplicity in $G_\mu(a_\mu|a)$. From the argument
in Case 2, we know that $x^{T}$ has a pole at most that of
$G_\mu(a_\mu|a)$.

If $\mu_{n+1-i}=\lambda_{n+1-i}$, then there will be a
corresponding factor $1+\beta a_k$ in $\Pi(a_\lambda)$ cancelling
that from $\Pi(a_\mu)$ ensuring that the multiplicity of $1+\beta
a_k$ in the denominator of $x^T\Pi(a_\mu)\Pi(a_\lambda)^{-1}$ is
not greater than that in $G_\lambda(a_\mu|a)$, as required.

Now we may suppose $\mu_{n+1-i}\neq \lambda_{n+1-i}$. We also must
have $G_\lambda(a_\mu|a)\neq 0$, so $\mu\supset\lambda$, and thus
$\mu_{n+1-i}>\lambda_{n+1-i}$.

Factors of $1+\beta a_k$ in the denominator of $x^T$ are in
one-to-one correspondence with occurrences of the entry $i$ in
$T$, so we only need to consider those $T$ with a maximal possible
occurrence of $i$ as given by Proposition \ref{prop}.

Consider the $\mu_{n+1-i}$-th column of our $\lambda$-tableau $T$,
and call it $C$.

As $(n+1-i,\mu_{n+1-i}+1)\not\in\mu$, $\mu'_{(\mu_{n+1-i}+1)}\leq
n-i$ and thus $T(1,\mu_{n+1-i}+1)>i$. Hence, there cannot be any
$i$'s to the right of $C$. Also, $n+1-\mu'_{\mu_{n+1-i}}\geq i$,
so since $T$ contains the maximal possible number of $i$'s, it
must contain an entry $i$ in the $\mu_{n+1-i}$-th column.

Let $j$ be the largest index for which $T(j,\mu_{n+1-i})$ contains
an entry less than $i+j$. Note $j$ must exist as there must exist
an $i$ in this column.

Pair the two tableaux $T$ and $T'$, identical in all respects
except that $T'$ contains an $i+j$ in the $j$-th row of $C$ and
$T$ does not. Note that these will both be semistandard since we
have the inequalities $T(j,\mu_{n+1-i}+1)\geq i+j$ (by Proposition
\ref{prop}) and $T(j+1,\mu_{n+1-i})>i+j$ (by maximality of $j$)
while $T(j,\mu_{n+1-i})$ already contains an entry less than
$i+j$.

We now calculate:

\begin{eqnarray*}
\beta^{|T|-|\lambda|}x^T+\beta^{|T'|-|\lambda|}x^{T'}&=&\beta^{|T|-|\lambda|}
x^T(1+\beta ( (a_\mu)_{n+1-(i+j)}) \oplus a_{i+j+c(j,\mu_{n+1-i})} )\\
&=&\beta^{|T|-|\lambda|} x^T \frac{1+\beta
a_{i+\mu_{n+1-i}}}{1+\beta a_{i+j-\mu_{n+1-(i+j)}}}
\end{eqnarray*}

By pairing our tableaux in this way and recovering an extra factor
$1+\beta a_k$ in the numerator, we see that the total pole for the
factor $1+\beta a_k$ is at most that of $G_\mu(a_\mu|a)$ as
required, so $c_\lambda^\mu$ cannot have any poles, since all
possible cases have now been considered.

Thus $c_\lambda^\mu$ is a polynomial in $\beta$ and $a$ (as we
already know it is a rational function of $\beta$ and $a$). We now
compute the degree of $c_\lambda^\mu$, considered as a polynomial
in $\beta$, and show by induction on $\mu$ that
$\deg_\beta{c_\lambda^\mu}\leq0$.

We use equation (\ref{pole}) and calculate the degree of each of
its constituent terms.

\begin{centre}
\begin{tabular}{|c|c|}
  \hline
  \hbox{term} & \mbox{degree} \\
  \hline
  $\beta^{|\mu|-|\lambda|}G_\mu(a_\mu|a)$ &
$|\mu|-|\lambda|-|\mu|=-|\lambda|$ \\[.2cm]
  $G_\lambda(a_\mu|a)$ & $\leq -|\lambda|$ \\[.2cm]
  $\Pi(a_\mu)\Pi(a_\lambda)^{-1}$ & 0 \\[.2cm]
  $\beta^{|\rho|-|\lambda|}c_\lambda^\rho G_\rho(a_\mu|a)$ & $\leq
|\rho|-|\lambda|+0-|\rho|=-|\lambda|$ \\[.2cm]
  \hline
\end{tabular}
\end{centre}

Here we use the fact that $x^T$ has degree $-|T|$ and the
inductive assumption for $\rho\subsetneq\mu$.

Now placing this into (\ref{pole}) we arrive at the inequality
$\deg{c_\lambda^\mu}\leq -|\lambda|+|\lambda|=0$ as required.

Being a polynomial in $\beta a_1,\beta a_2,\ldots $ of degree at
most zero in $\beta$, $c_\lambda^\mu$ must be constant, that is
independent of $\beta$ and $a$. Thus we can calculate the values
of $c^\mu_\lambda$ by specialisation to the ordinary Grothendieck
polynomials with $a=0$. From Proposition \ref{buch}, we know the
value of $c_\lambda^\mu(\beta,0)$ and thus,
$$c_\lambda^\mu(\beta,a)=c_\lambda^\mu(\beta,0)=\begin{cases}
        1,
                    &\text{if\ } \lambda\tangle\mu,\\
            0,&\text{otherwise}.
       \end{cases}$$ so Proposition \ref{main1prop} is
proven, as required.

\subsection{The recurrence relation}

Suppose $\nu$ is a partition of length at most $n$, $\mu$ is a
partition, $\theta$ is a skew diagram and $P(x)$ is a fixed
symmetric polynomial in $x_1,x_2\ldots x_n$ with coefficients in
$\F(\beta,a)$. Then define the coefficients
$g^\nu_{\mu}=g^\nu_{\mu}(P)\in \F(\beta,a)$ by the formula
\begin{equation} \label{coef}
P(x)G_\mu(x|a)=\sum_\nu g^\nu_{\mu} G_\nu(x|a).
\end{equation} Theorem \ref{basis} ensures that these coefficients
are well defined.

In the important special case where $P(x)=G_\th(x|b)$ with
$b=(b_i)_{i\in\Z}$ a second doubly infinite sequence of variables,
we use the notation $g_{\th\mu}^\nu=g_{\th\mu}^\nu(a,b)$ for
$g_\mu^\nu(G_\th(x|b))$.

\begin{proposition} The coefficients $g_{\mu}^\nu$ satisfy the following
recurrence:

\begin{equation}\label{recurrence}
g_{\mu}^{\nu}=\frac{1}{\Pi(a_{\nu})-\Pi(a_{\mu})}
\Big(\sum_{\mu\tangle^*\lambda}\Pi(a_{\mu})\beta^{|\lambda/\mu|}g_{\lambda}^{\nu}
-\sum_{\eta\tangle^*\nu}\Pi(a_{\eta})\beta^{|\nu/\eta|}g_{\mu}^{\eta}
\Big),
\end{equation} with boundary conditions \begin{eqnarray*}(i)&&
g^\nu_{\mu}=0\quad\hbox{unless}\quad\mu\subset\nu, \\ (ii)&&
g^\lambda_{\lambda}=P(a_\lambda).\end{eqnarray*}\end{proposition}

This is indeed a recurrence, for it enables the coefficients
$g^\nu_{\mu}$ to be computed recursively by induction on
$|\nu/\mu|$.

\begin{proof} Applying Proposition \ref{main1prop} to
\[
P(x)G_\mu(x|a)\Pi(x)=\sum_\eta g^\eta_{\mu} G_\eta(x|a)\Pi(x)
\] yields the following:
\[
P(x)\Pi(a_\mu)\sum_{\mu\tangle\lambda}\beta^{|\lambda/\mu|}G_\lambda(x|a)=
\sum_\eta
g_{\mu}^\eta\Pi(a_\eta)\sum_{\eta\tangle\nu}\beta^{|\nu/\eta|}G_\nu(x|a).
\]
If we now combine this with (\ref{coef}) we obtain the identity
\[
\Pi(a_\mu)\sum_{\mu\rightrightarrows\lambda}\sum_\nu\beta^{|\lambda/\mu|}g^\nu_{\lambda}
G_\nu(x|a)=\sum_\eta
g_{\mu}^\eta\Pi(a_\eta)\sum_{\eta\tangle\nu}\beta^{|\nu/\eta|}G_\nu(x|a).
\]
We now use the fact that the factorial Grothendieck polynomials
$G_\lambda(x|a)$ form a basis to equate the coefficients of
$G_\nu(x|a)$, giving
$$\Pi(a_\mu)\sum_{\mu\tangle\lambda}\beta^{|\lambda/\mu|}g^\nu_{\lambda}
=\sum_{\eta\tangle\nu}\Pi(a_\eta)\beta^{|\nu/\eta|}g^\eta_{\mu}$$
which rearranges to the quoted form of the recurrence.

For the boundary conditions, suppose that $\rho$ is a minimal
partition with respect to containment order such that
$g_{\mu}^\rho\neq 0$. Substituting $x=a_\rho$ in (\ref{coef})
gives $G_\mu(a_\rho|a)P(a_\rho)=g_{\mu}^\rho G_\rho(a_\rho|a)$. If
$\mu\not\subset\rho$ then from the vanishing theorem, we get $(i)$
$g_{\mu}^\rho=0$, so now we may deal with the $\mu=\rho$ case
which gives $(ii)$ $g_{\rho}^\rho=P(a_\rho)$ as required.
\end{proof}

We now give a general solution to the above recurrence.

For a partition $\lambda$, introduce the notation $\Pi(\lambda)$
to represent $\Pi(a_\lambda)$.

\begin{proposition}

The general solution to the recurrence (\ref{recurrence}) is
\[
g_\mu^\nu=\beta^{|\nu/\mu|}\sum_R
\Pi(\rho_0)\Pi(\rho_1)\ldots\Pi(\rho_{l-1})\sum_{k=0}^l
P(a_{\rho_k})\prod_{\substack{i=0 \\ i\neq k}}^l
\frac{1}{\Pi(\rho_k)-\Pi(\rho_i)}
\]
where the sum is over all sequences
$$R:\mu=\rho_0\tangle^*\rho_1\tangle^*\cdots\tangle^*\rho_{l-1}\tangle^*\rho_l=\nu.$$ \end{proposition}

\begin{proof} We need to show that this proposed solution
satisfies both the recurrence relation and the boundary
conditions.

That this proposed solution satisfies the boundary conditions is
immediate, for if $\mu\not\subset\nu$ there is no such sequence
$R$ while if $\mu=\nu$ there is exactly one such sequence, the
trivial sequence with $l=0$.

By induction on $|\nu/\mu|$, we get

\begin{eqnarray*}
g_\mu^\nu&=&\frac{1}{\Pi(\nu)-\Pi(\mu)}
\Big(\sum_{\mu\tangle^*\lambda}\Pi(a_{\mu})\beta^{|\lambda/\mu|}g_{\lambda}^{\nu}
-\sum_{\eta\tangle^*\nu}\Pi(a_{\eta})\beta^{|\nu/\eta|}g_{\mu}^{\eta}
\Big) \\
&=&\frac{\beta^{|\nu/\mu|}}{\Pi(\nu)-\Pi(\mu)} \Big(\sum_R
\prod_{j=0}^{l-1}\Pi(\r{j}) \sum_{k=1}^{l}
P(a_{\rho_k})\prod_{\substack{i=1 \\ i\neq k}}^l
\frac{1}{\Pi(\rho_k)-\Pi(\rho_i)}\\ &&\qquad\qquad{}-\sum_R
\prod_{j=0}^{l-1}\Pi(\r{j})\sum_{k=0}^{l-1}
P(a_{\rho_k})\prod_{\substack{i=0 \\ i\neq k}}^{l-1}
\frac{1}{\Pi(\rho_k)-\Pi(\rho_i)} \Big) \\
&=&\beta^{|\nu/\mu|} \sum_R \prod_{j=0}^{l-1}\Pi(\r{j})
\sum_{k=0}^{l}
\frac{P(a_{\r{k}})[(\Pi(\rho_k)-\Pi(\rho_0))-(\Pi(\rho_k)-\Pi(\rho_l))]}{(\Pi(\nu)-\Pi(\mu))\prod_{\substack{i=0
\\
i\neq k}}^l \Pi(\rho_k)-\Pi(\rho_i)} \\
&=&\beta^{|\nu/\mu|}\sum_R
\Pi(\rho_0)\Pi(\rho_1)\ldots\Pi(\rho_{l-1})\sum_{k=0}^l
P(a_{\rho_k})\prod_{\substack{i=0 \\ i\neq k}}^l
\frac{1}{\Pi(\rho_k)-\Pi(\rho_i)}
\end{eqnarray*} as required. \end{proof}

\section{Calculation of the Coefficients}

The general solution to the recurrence appears inadequate, in that
it is hard to specialise to the case of ordinary Grothendieck
polynomials by setting $a=0$, and nor does it clearly reflect the
stringent conditions we have imposed on denominators in Lemma
\ref{polecontrol}. So now we turn specifically to the case
$P(x)=G_\theta(x|b)$ and provide an alternative description of the
coefficients $g_{\th\mu}^\nu$ with a view to specialising to the
ordinary Grothendieck polynomials.

\subsection{Solution where all boxes of $\theta$ are in different
columns}

We now provide a solution to the recurrence in the case where
$\theta$ does not contain two boxes in the same column. In order
to state this result however, we first need to define some more
combinatorial objects.

Consider a sequence of partitions
\begin{equation}\label{picseq}
R:\mu=\rho^{(0)}\o{1}\rho^{(1)}\o{2}\cdots\o{l-1}\rho^{(l-1)}\o{l}\rho^{(l)}=\nu.
\end{equation}

Say a semistandard set-valued $\theta$ tableau $T$ is {\it
related\/} to $R$ if $T$ contains distinguished entries
$r_1,r_2,\ldots,r_l$ in cells $\alpha_1,\alpha_2,\ldots,\alpha_l$
respectively with $(r_1,\alpha_1)\prec
(r_2,\alpha_2)\prec\ldots\prec (r_l,\alpha_l)$ where $\prec$ is
the ordering defined in Section \ref{tableaux}. We distinguish
these entries by placing a bar over them.

If $\r{0},\r{1},\ldots,\r{l}$ are partitions,
$r_1,r_2\ldots,r_l\in[n]$, and $T$ is a $\theta$-tableau with
distinguished entries
$(r_1,\alpha_1)\prec\cdots\prec(r_l,\alpha_l)$, then we define the
function
$$F_T\Big(\r{l}\mathop{|}_{r_l}\cdots\m{2}\r{1}\m{1}\r{0}\Big)$$
to equal the product
$$\prod_{\substack{\alpha \in\theta \\ r\in T(\alpha)\\ {\rm
unbarred}}}((a_{\rho(r)})_r\oplus
b_{r+c(\alpha)})\cdot\prod_{i=1}^l(1+\beta
(a_{\r{i-1}})_{r_i})(1+\beta b_{r_i+c(\alpha_i)}).$$ where
$\rho(r)=\r{i}$ if $r_i\prec r\prec r_{i+1}$. In the important
case where the $\r{i}$ and $r_i$ form a sequence $R$ of the form
(\ref{picseq}), and $T$ is a semistandard $\theta$-tableau related
to $R$, we denote this product by $w(T)$.

\begin{theorem}\label{bsolnthm} For $P(x)=G_{\theta}(x|b)$, if $\theta$
does not contain two boxes in the same column, then
\begin{equation} \label{bsoln} g_{\theta\mu}^{\nu}=\sum_{(R,T)}
\beta^{|T|-|\theta|}w(T).
\end{equation} where the sum is over all $\theta$-tableaux $T$ which are
related to a sequence $R$ of the form (\ref{picseq}).
\end{theorem}

\begin{proof} For $\nu\not\supset\mu$, no such sequences $R$ exist so
(\ref{bsoln}) agrees with $g_{\theta\mu}^\nu=0$ as required.

For $\nu=\mu$, there is one such sequence $R$, $\r{0}=\r{l}=\mu$,
so no barred entries can exist in $T$. The set of tableaux summed
over is now exactly the same as the set summed over in the
definition of $G_\theta(x|a)$, and we thus notice that
$\sum_{(R,T)}\beta^{|T|-|\theta|}w(T)=G_\theta(a_\lambda|b)$ while
the boundary conditions of the recurrence give
$g_{\theta\mu}^\nu=G_\theta(a_\lambda|b)$, agreeing with
(\ref{bsoln}) as required.

Now we need to show that our proposed solution satisfies the
recurrence. So we suppose (\ref{bsoln}) holds and we have to show
that this implies (\ref{recurrence}) holds.

Let $m$ be a non-negative integer. Let $l=|\nu/\mu|$. We now form
a set $\mathcal{T}_m$ of triples $(k,R,T)$ as follows.

$k$ is an integer from $m$ to $l$ inclusive. $R$ is a sequence
$$R:\mu=\rho^{(0)}\o{1}\ldots\o{k-m}\rho^{(k-m)}\tangle\rho^{(k)}\o{k+1}\ldots\o{l}\rho^{(l)}=\nu.$$

$T$ is a semistandard set-valued $\theta$-tableau $T$ with entries
from $[n]$ such that $T$ contains distinguished entries
$r_1,r_2,\ldots,r_{k-m},r_{k+1},\ldots,r_l$ with $r_1\prec
r_2\prec\ldots\prec r_{k-m}\prec r_{k+1}\prec\ldots\prec r_l$.
These entries are distinguished by placing a bar over them.

$\mathcal{T}_m=\mathcal{T}_m(\theta,\mu,\nu)$ is defined to be the
set of all such triples $(k,R,T)$ as defined above.

We define two weights on such a triple $(k,R,T)$, a positive and a
negative weight, by

$$w^+(k,R,T)=\beta^{|T|-|\theta|}F_T\Big(\r{l}\m{l}\cdots\m{k+1}\r{k}\m{k-m}\cdots\m{1}\r{0}\Big),$$
$$w^-(k,R,T)=\beta^{|T|-|\theta|}F_T\Big(\r{l}\m{l}\cdots\m{k+1}\r{k-m}\m{k-m}\cdots\m{1}\r{0}\Big)
\frac{1+\beta (a_{\r{k}})_{r_{k+1}}}{1+\beta
(a_{\r{k-m}})_{r_{k+1}}}.$$ The extra factor in the definition of
$w^-$ is to \lq correct' the contribution to the product provided
by the barred $r_{k+1}$.

Let
$$S_m=\sum_{(k,R,T)\in\mathcal{T}_m}w^+(k,R,T)\Pi(\r{k-m})-\sum_{(k,R,T)\in\mathcal{T}_m}w^-(k,R,T)\Pi(\r{k-m}).$$

Say that $(k,R,T)$ is a positive $\varepsilon$-triple if
$\r{k}/\r{k-m-1}$ contains the shape
$\setlength{\unitlength}{0.5em}
\begin{picture}(2,2)
\put(0,1){\line(1,0){2}} \put(0,2){\line(1,0){2}}
\put(1,0){\line(1,0){1}} \put(2,2){\line(0,-1){2}}
\put(0,1){\line(0,1){1}} \put(1,2){\line(0,-1){2}}
 \end{picture}
\setlength{\unitlength}{1pt}$. Say that $(k,R,T)$ is a negative
$\varepsilon$-triple if $\r{k+1}/\r{k-m}$ contains the shape
$\setlength{\unitlength}{0.5em}
\begin{picture}(2,2)
\put(0,1){\line(1,0){2}} \put(0,2){\line(1,0){1}}
\put(0,0){\line(1,0){2}} \put(2,1){\line(0,-1){1}}
\put(0,0){\line(0,1){2}} \put(1,2){\line(0,-1){2}}
 \end{picture}
\setlength{\unitlength}{1pt}$. Let $\T_m^+=\{ (k,R,T)\in\T_m \mid
(k,R,T) \text{ is a positive $\epsilon$-triple} \}$ and similarly
$T_m^-=\{ (k,R,T)\in\T_m \mid (k,R,T) \text{ is a negative
$\epsilon$-triple} \}$. Define
$$\varepsilon_m=\sum_{(k,R,T)\in\T^+_m}w^+(k,R,T)\Pi(\r{k-m})-
\sum_{(k,R,T)\in\T^-_m}w^-(k,R,T)\Pi(\r{k-m}).$$

$S_0=0$, as $w^+$ and $w^-$ are identical functions when $m=0$.
$S_{l+1}=0$, as there are no sequences with $m=l+1$. Similarly
$\varepsilon_0=\varepsilon_{m+1}=0$. Now if we temporarily assume
Proposition \ref{smprop} below, we can obtain the equation
$$\sum_{\mu\tangle\lambda}\beta^{|\lambda/\mu|}g^\nu_{\theta\lambda}\Pi(\mu)
=\sum_{\eta\tangle\nu}\beta^{|\nu/\eta|}g^\eta_{\th\mu}\Pi(\eta)$$
which is equivalent to the recurrence, so we are done. \end{proof}

So the proof of Theorem \ref{bsolnthm} follows immediately from
the proof of the following proposition.

\begin{proposition} \label{smprop}
\begin{equation} \label{sm}
S_m-\varepsilon_m=\beta(S_{m+1}-\varepsilon_{m+1})+\sum_{\substack{\mu\tangle\lambda\\
|\lambda/\mu|=m}}g^\nu_{\th\lambda}
\Pi(\mu)-\sum_{\substack{\eta\tangle\nu\\
|\nu/\eta|=m}}g^\eta_{\th\mu} \Pi(\eta)
\end{equation} \end{proposition}


\begin{proof}
The positive terms in $S_m$ with $k=m$ and the negative terms in
$S_m$ with $k=l$ give exactly $$\sum_{\substack{\mu\tangle\lambda\\
|\lambda/\mu|=m}}g^\nu_{\th\lambda}
\Pi(\mu)-\sum_{\substack{\eta\tangle\nu\\
|\nu/\eta|=m}}g^\eta_{\th\mu} \Pi(\eta).$$

So now we consider positive terms in $S_m$ with $k>m$, and
negative terms with $k<l$. For positive terms, we consider
$\Theta=\r{k}/\r{k-m-1}$ while for negative terms we consider
$\Theta=\r{k+1}/\r{k-m}$. We have two separate cases to consider,
according to the shape of $\Theta$.

{\bf Case 1:} $\Theta$ contains two boxes in the same row.

Consider a positive term $w^+(k,R,T)$ covered by this case. Define
$(k',R',T')$ as follows: Set $k'=k-1$. Construct $R'$ from $R$ by
replacing the subsequence $\r{k-m-1}\o{k-m}\r{k-m}\tangle\r{k}$ by
$\r{k-m-1}\tangle\rho'\o{k-m}\r{k}$ (there exists a unique such
partition $\rho'$). Set $T'=T$.

Then $w^+(k,R,T)\Pi(\r{k-m})=w^-(k',R',T')\Pi(\r{k'-m})$. For the
only factor differing in $w^+(k,R,T)$ and $w^-(k',R',T')$ is that
due to the barred $r_{k-m}$. In $w^+(k,R,T)$, this entry
contributes the factor $(1+\beta(a_{\r{k-m-1}})_{r_{k-m}})(1+\beta
b_{r_{k-m}+c(\alpha)})$ while in $w^-(k',R',T')$, this entry
contributes $(1+\beta(a_{\rho'})_{r_{k-m}})(1+\beta
b_{r_{k-m}+c(\alpha)})$, precisely countering the difference in
the factors $\Pi(\r{k-m})$ and $\Pi(\r{k'-m})$.

This map $(k,R,T)\mapsto(k',R',T')$ has a similar inverse, hence
is bijective, so we have shown that all terms of $S_m$ which are
covered by this case cancel each other to give no net
contribution.

{\bf Case 2:} All boxes of $\Theta$ are in different rows and
columns.

Given such a positive triple $(k_0,R_0,T_0)$, consider all such
triples $(k,R,T)$ with $k=k_0$, $T=T_0$ and $R=R_0$ except for
$\r{k-m}$ (so there are $m+1$ such triples). Also consider all
negative triples $(k',R',T')$ with $k'=k-1$, $T'=T_0$ and $R=R_0$
except that the subsequence $\r{k-m-1}\to\r{k-m}\tangle\r{k}$ is
replaced by $\r{k-m-1}\tangle\rho'\to\r{k}$ for one of the $m+1$
possibilities for $\rho'$.

Let the row numbers of $\r{k}/\r{k-m-1}$ be
$s_1,s_2,\ldots,s_{m+1}$. Let $y_i=(a_{\r{k}})_{i}$,
$z_i=(a_{\r{k-m-1}})_{i}$.

Then these $2m+2$ triples together contribute the following to the
sum $S_m$:

\[
\beta^{|T|-|\theta|}\sum_{j=1}^{m+1}F_T\Big(
\cdots\mathop{|}\r{k}\mathop{|}_{s_j}\r{k-m-1}\mathop{|}\cdots
\Big)\Big( \frac{1+\beta y_{s_j}}{1+\beta z_{s_{j}}}
\Big)\Pi(\r{k-m-1})
\]

Between $\overline{r_{k+1}}$ and $\overline{r_{k-m-1}}$, suppose
the $s_j$'s (all possible $j$) occur in order $t_1\prec
t_2\prec\cdots\prec t_p$ and suppose $t_i$ lies in cell
$\alpha_i$.

We only need to consider the entries $t_1,\ldots,t_p$ in $T$, for
all other entries, along with $\Pi(\r{k-m-1})$ contribute a common
factor. After taking out that very common factor, and noticing
that the relevant barred $s_j$ contributes a factor $(1+\beta
z_{s_j})(1+\beta b_{s_j+c(\alpha_j)})$, we get
\begin{eqnarray*}
&\ &\sum_{i=1}^p \beta(y_{t_i}-z_{t_i})(1+\beta
b_{t_i+c(\alpha_i)}) \prod_{j=1}^{i-1}y_{t_j}\oplus
b_{t_j+c(\alpha_j)}
\prod_{j=i+1}^{p}z_{t_j}\oplus b_{t_j+c(\alpha_j)} \\
&=&\beta\sum_{i=1}^p (y_{t_i}\oplus
b_{t_i+c(\alpha_i)}-z_{t_i}\oplus
b_{t_i+c(\alpha_i)})\prod_{j=1}^{i-1}y_{t_j}\oplus
b_{t_j+c(\alpha_j)} \prod_{j=i+1}^{p}z_{t_j}\oplus
b_{t_j+c(\alpha_j)}.
\end{eqnarray*}

This is a telescoping sum and equals
$$\beta\Big( \prod_{j=1}^{p}y_{t_j}\oplus
b_{t_j+c(\alpha_j)}-\prod_{j=1}^{p}z_{t_j}\oplus
b_{t_j+c(\alpha_j)} \Big).$$

Now if we replace the common factor, we obtain
$$\beta(w^+(k,R^*,T)-w^-(k,R^*,T))\Pi(\r{k-(m+1)})$$ where $R^*$ is
the sequence obtained by replacing
$\r{k-m-1}\to\r{k-m}\tangle\r{k}$ by $\r{k-(m+1)}\tangle\r{k}$.

Hence, when considering the contribution of all terms of $S_m$
covered by this case, they add up to give exactly $\beta S_{m+1}$.

{\bf Case 3:} $\Theta$ contains two boxes in the same column, but
does not contain two boxes in the same row.

Let $i$ and $i+1$ be the row numbers of the two boxes of $\Theta$
which are in the same column. We will underline the marked $i$ and
$i+1$ in $T$ which come from $\Theta$ for increased clarity.

Tableaux containing the following cannot occur as they cannot
arise from a sequence of partitions.

\[
\setlength{\unitlength}{17pt}
\begin{picture}(3,1) \put(0,0){\line(0,1){1}}
\put(0,0){\line(1,0){3}} \put(1,0){\line(0,1){1}}
\put(2,0){\line(0,1){1}} \put(0,1){\line(1,0){3}}
\put(3,0){\line(0,1){1}}
\put(0.4,0.25){$\underline{\overline{i}}$}
\put(2.4,0.25){$\overline{i}$} \put(1.1,0.25){$\cdots$}
\end{picture}
\setlength{\unitlength}{1pt} \qquad \setlength{\unitlength}{17pt}
\begin{picture}(3,1) \put(0,0){\line(0,1){1}}
\put(0,0){\line(1,0){3}} \put(1,0){\line(0,1){1}}
\put(2,0){\line(0,1){1}} \put(3,0){\line(0,1){1}}
\put(0,1){\line(1,0){3}}
\put(0.067,0.25){$\overline{i\!\!+\!\!1}$}
\put(2.067,0.25){$\underline{\overline{i\!\!+\!\!1}}$}
\put(1.1,0.25){$\cdots$}
\end{picture}
\setlength{\unitlength}{1pt}\qquad
\]

If $\underline{\overline{i}}$ appears to the left of
$\overline{i}$ in the same row, then by examining $R$, there would
need to be a marked $i+1$ between these two entries in reverse
column order. But since $\theta$ does not contain two boxes in the
same column, this cannot happen for the tableau to be
semistandard. The case of $\underline{\overline{i+1}}$ to the
right of $\overline{i+1}$ in the same row proceeds similarly.

Wherever possible, we match up our tableaux as follows:

Given an $\underline{\overline{i+1}}$ in a negative term, denote
by $L$ the box with this entry and all consecutive boxes to its
left which contain an $i+1$. Consider the two tableaux $T_1$ and
$T_2$, the first with an unbarred $i$ in the first box of $L$ and
the second without this $i$. We match these up with $T_3$ and
$T_4$ which are obtained by changing all $i$'s and $i+1$'s in $L$
such that $\underline{\overline{i}}$ is in the leftmost box of
$L$, an unbarred $i$ is in each other box of $L$ and $T_3$
contains an unbarred $i+1$ in the rightmost box of $L$ while $T_4$
does not. For example (i=2):
\[
T_1\hspace{22pt}+\hspace{22pt}
T_2\hspace{20pt}\quad\longleftrightarrow \hspace{20pt}\quad
T_3\hspace{22pt}+\hspace{22pt} T_4
\]
$$\setlength{\unitlength}{1.6em}
\parbox[c]{7\unitlength}{\begin{picture}(7,1)
\multiput(0,0)(1,0){8}{\line(0,1){1}}
\multiput(0,0)(0,1){2}{\line(1,0){3}}
\multiput(4,0)(0,1){2}{\line(1,0){3}} \put(0.2,0.3){$23$}
\put(1.36,0.3){$3$} \put(2.36,0.3){$\ou{3}$} \put(3.36,0.3){$+$}
\put(4.36,0.3){$3$} \put(5.36,0.3){$3$} \put(6.36,0.3){$\ou{3}$}
\end{picture}}
\quad\longleftrightarrow\quad
\parbox[c]{7\unitlength}{\begin{picture}(7,1)
\multiput(0,0)(1,0){8}{\line(0,1){1}}
\multiput(0,0)(0,1){2}{\line(1,0){3}}
\multiput(4,0)(0,1){2}{\line(1,0){3}} \put(0.36,0.3){$\ou{2}$}
\put(1.36,0.3){$2$} \put(2.2,0.3){$23$} \put(3.36,0.3){$+$}
\put(4.36,0.3){$\ou{2}$} \put(5.36,0.3){$2$} \put(6.36,0.3){$2$}
\end{picture}}
\setlength{\unitlength}{1pt}
$$

Now we shall show that under this identification, the
corresponding terms of $S_m$ cancel, that is
\begin{equation}\label{theend}
w^-(k-1,R',T_1)+w^-(k-1,R',T_2)=w^+(k,R,T_3)+w^+(k,R,T_4)
\end{equation} where $R'$ is the sequence obtained from $R$ by
replacing $\r{k-m-1}\to\r{k-m}\tangle\r{k}$ with
$\r{k-m-1}\tangle\rho^{'} \to\r{k}$.

This is because the two sides of the equation have common factors
from their common entries, as well as from the unmarked $i$'s in
the positive terms and the unmarked $i+1$'s in the negative terms.
Apart from these common factors, the positive terms have, upon
combination, the extra factors $(1+\beta v)(1+\beta
b_{i+c(\alpha_l)})$ from $\underline{\overline{i}}$, $(1+\beta
w)(1+\beta b_{i+1+c(\alpha_r)})$ from $i+1$ and $1+\beta u$ from
$\Pi(\r{k-m})$, while the negative terms have the extra factors
$(1+\beta w)(1+\beta b_{i+1+c(\alpha_r)})$ from
$\underline{\overline{i+1}}$, $(1+\beta u)(1+\beta
b_{i+c(\alpha_l)})$ from $i$ and $1+\beta v$ from
$\Pi(\r{k-m-1})$. Here $u=(a_{\r{k}})_i$, $v=(a_{\r{k}})_{i+1}$,
$w=(a_{\r{k-m}})_{i+1}$ and $\alpha_l$ and $\alpha_r$ are
respectively the leftmost and rightmost boxes of $L$.

Taking into account those tableaux which we have already shown to
give zero contribution to the sum, we find that the only remaining
tableaux $T$ for which the above identification of positive and
negative terms cannot be made, is where a barred $i$ is in the
leftmost cell of $L$ ($L$ as defined above) for a negative term
and vice versa for a positive term. We shall focus our attention
on the positive terms for which this happens, as the negative case
proceeds similarly.

Then the following must occur as a subsequence of $R$:

$$\r{k-m-2}\stackrel{i+1}\longrightarrow\r{k-m-1}\stackrel{i}\longrightarrow\r{k-m}\tangle\r{k}.$$

By replacing this with
$\r{k-m-2}\tangle\rho'\stackrel{i+1}\longrightarrow\r{k}$, we form
another sequence $R'$.

We map $T$ with $i$ barred in the leftmost box of $L$ to $T'$ and
$T''$ with the barred $i$ removed in the first case and unbarred
in the second case. All boxes of $L$ in $T'$ and $T''$ contain an
$i+1$ with the rightmost of these barred. For example:

\[
\setlength{\unitlength}{1.6em}
\parbox[c]{3\unitlength}{\begin{picture}(3,1)
\multiput(0,0)(1,0){4}{\line(0,1){1}}
\multiput(0,0)(0,1){2}{\line(1,0){3}} \put(0.36,0.3){$\ou{2}$}
\put(1.36,0.3){$2$} \put(2.2,0.3){$2\overline{3}$}
\end{picture}}
\mapsto
\parbox[c]{7\unitlength}{\begin{picture}(7,1)
\multiput(0,0)(1,0){8}{\line(0,1){1}}
\multiput(0,0)(0,1){2}{\line(1,0){3}}
\multiput(4,0)(0,1){2}{\line(1,0){3}} \put(0.36,0.3){$3$}
\put(1.36,0.3){$3$} \put(2.36,0.3){$\overline{3}$}
\put(3.36,0.3){$+$} \put(4.2,0.3){$23$} \put(5.36,0.3){$3$}
\put(6.36,0.3){$\overline{3}$}
\end{picture}}
\setlength{\unitlength}{1pt}
\]

\begin{lemma} \label{claim1}
\[
w^+(k,R,T)\Pi(\r{k-m})=\beta(w^-(k-1,R',T')+w^-(k-1,R',T''))\Pi(\r{k-m-2}).
\]
\end{lemma}

\begin{proof} Define $u,v,w,x$ by the following, where the pronumeral in column
$\rho$ and row $j$ represents $(a_\rho)_j$.

\[
\begin{array}{c|c|c|c|c}
  \  & \r{k-m-2} & \r{k-m-1},\rho' & \r{k-m} & \r{k} \\ \hline
  i & v & v & u & u \\
  i+1 & x & w & w & v \\
\end{array}
\]

Each side of (\ref{claim1}) has common factors due to common
entries and from the unbarred $i$'s and $i+1$'s in $L$ in the
positive and negative terms respectively. Aside from these common
factors, the left hand side of (\ref{claim1}) has factors
$(1+\beta v)(1+\beta b_{i+c(\alpha_l)})$ from
$\underline{\overline{i}}$, $(1+\beta x)(1+\beta
b_{i+1+c(\alpha_r)})$ from $\overline{i+1}$ and $(1+\beta
u)(1+\beta w)$ from $\Pi(\r{k-m})$, while the right hand side has
factors $(1+\beta w)(1+\beta b_{i+1+c(\alpha_r)})$ from
$\overline{i+1}$, $(1+\beta u)(1+\beta b_{i+c(\alpha_l)})$ from
$i$ and $(1+\beta v)(1+\beta x)$ from $\Pi(\r{k-m-2})$. Again
$\alpha_l$ and $\alpha_r$ denote the leftmost and rightmost cells
of $L$. The extra factor of $\beta$ is due to there being one more
entry in $T$ than in $T'$.
\end{proof}

We note that the inverse to this map can always be created, for
given a tableau related to a sequence of the form $R'$, analysis
of the sequence of partitions shows that a barred $i+1$ can never
be to the left of the relevant barred $i+1$, ensuring an
uninterrupted string of unbarred $i+1$'s to the left of the
relevant barred $i+1$ which enables the inverse to be easily
constructed. Thus positive terms covered by Case 3 in $S_m$ give
exactly $\beta$ times the negative terms in $\varepsilon_{m+1}$.
Similarly, we see that the negative terms covered by Case 3 in
$S_m$ give $\beta$ times the positive terms in
$\varepsilon_{m+1}$.

{\bf Case 4:} The only remaining triples are
$\varepsilon$-triples, so their contribution to $S_m$ is exactly
$\varepsilon_m$.

Hence we have proven (\ref{sm}), as desired. \end{proof}

\subsection{Partial solution in the general case}

We now consider the case where $b=0$, with a view to turning our
attention to the ordinary Grothendieck polynomials. We return to
the situation where $\theta$ is an arbitrary skew diagram. In this
case, we provide a partial solution to the recurrence relation. We
shall carry over notation used in the case of arbitrary $b$, just
noting that the variables $b_i$ are all to be set equal to zero.

We shall also set the following variables $a_i$ equal to zero: If
there exists a sequence $R$ of the form (\ref{picseq}), and a $k$
for which $\rho^{(k+1)}/\rho^{(k-1)}$ consists of two boxes in the
same column, then set $(a_{\rho^{(k-1)}})_{r_k}=0.$ Let us pause
and note that this is equivalent to $(a_{\r{k+1}})_{r_{k+1}}=0$.

\begin{theorem} If the appropriate variables are all set to zero as
described above, then we have
\begin{equation} \label{zsoln} g_{\th \mu}^{\nu}(a,0)=\sum_{(R,T)}
\beta^{|T|-|\theta|}w(T).\end{equation} where again, the sum is
over all $\theta$-tableaux $T$ which are related to a sequence $R$
of the form (\ref{picseq}).
\end{theorem}

\begin{proof} As in the proof of Theorem \ref{bsolnthm}, this proposed solution satisfies the
boundary conditions of the recurrence. So now we suppose that
$\mu\subsetneq\nu$. Despite setting some of the variables $a_i$
equal to zero, we still have $\Pi(\nu)\neq\Pi(\mu)$, so we are
able to calculate the coefficients using the recurrence
(\ref{recurrence}). So it now suffices to show that our proposed
solution satisfies the recurrence. So we suppose that
(\ref{zsoln}) holds and use this to show that (\ref{recurrence})
holds.

Again as in the proof of Theorem \ref{bsolnthm}, the proof reduces
to the proof of the following proposition. \end{proof}

\begin{proposition} \label{zemprop}
\begin{equation} \label{zem}
S_m-\varepsilon_m=\beta(S_{m+1}-\varepsilon_{m+1})+\sum_{\substack{\mu\tangle\lambda\\
|\lambda/\mu|=m}}g^\nu_{\th\lambda}
\Pi(\mu)-\sum_{\substack{\eta\tangle\nu\\
|\nu/\eta|=m}}g^\eta_{\th\mu} \Pi(\eta)
\end{equation} \end{proposition}

\begin{proof} As per the proof of Proposition \ref{smprop}, in $S_m$,
the positive terms with $k=m$ and the negative terms with $k=l$
give exactly
$$\sum_{\substack{\mu\tangle\lambda\\
|\lambda/\mu|=m}}g^\nu_{\th\lambda}
\Pi(\mu)-\sum_{\substack{\eta\tangle\nu\\
|\nu/\eta|=m}}g^\eta_{\th\mu} \Pi(\eta)$$

For the remaining terms in $S_m$, we split them up into three
cases according to the shape of $\Theta$.

{\bf Case 1:} $\Theta$ contains two boxes in the same row, and
does not contain two boxes in the same column.

This case is the same as in the previous solution, all such
$\Theta$ combined contribute zero to the sum.

{\bf Case 2:} $\Theta$ contains all boxes in different rows and
columns.

Again this case is the same as in the previous solution,
contributing $\beta S_{m+1}$ to the sum.

{\bf Case 3:} $\Theta$ contains two boxes in the same column, but
does not contain two boxes in the same row.

Let $i$ and $i+1$ be the row numbers of the two boxes of $\Theta$
which are in the same column. We will underline the marked $i$ and
$i+1$ in $T$ which come from $\Theta$ for increased clarity.

Tableaux containing the following have zero weight, due to setting
variables equal to zero, so their contribution can be neglected:

\[
\setlength{\unitlength}{17pt}
\begin{picture}(1,2)
\put(0,0){\line(1,0){1}} \put(0,0){\line(0,1){2}}
\put(0,1){\line(1,0){1}} \put(0,2){\line(1,0){1}}
\put(1,0){\line(0,1){2}}
\put(0.4,1.25){$\underline{\overline{i}}$}
\put(0.1,0.25){$i\!\!+\!\!1$}
\end{picture}
\setlength{\unitlength}{1pt}\qquad\setlength{\unitlength}{17pt}
\begin{picture}(2,1)
\put(0,0.5){\line(0,1){1}} \put(0,0.5){\line(1,0){2}}
\put(1,0.5){\line(0,1){1}} \put(2,0.5){\line(0,1){1}}
\put(0,1.5){\line(1,0){2}}
\put(0.4,0.75){$\underline{\overline{i}}$} \put(1.4,0.75){$i$}
\end{picture}
\setlength{\unitlength}{1pt}\qquad\setlength{\unitlength}{17pt}
\begin{picture}(1,2) \put(0,0){\line(1,0){1}}
\put(0,0){\line(0,1){2}} \put(0,1){\line(1,0){1}}
\put(0,2){\line(1,0){1}} \put(1,0){\line(0,1){2}}
\put(0.4,1.25){$i$}
\put(0.05,0.25){$\underline{\overline{i\!\!+\!\!1}}$}
\end{picture}
\setlength{\unitlength}{1pt}\qquad \setlength{\unitlength}{17pt}
\begin{picture}(2,1) \put(0,0.5){\line(0,1){1}}
\put(0,0.5){\line(1,0){2}} \put(1,0.5){\line(0,1){1}}
\put(2,0.5){\line(0,1){1}} \put(0,1.5){\line(1,0){2}}
\put(0.07,0.75){$i\!\!+\!\!1$}
\put(1.07,0.75){$\underline{\overline{i\!\!+\!\!1}}$}
\end{picture}
\setlength{\unitlength}{1pt}
\]

This is because, in the first case, the $i+1$ contributes
$(a_{\r{k}})_{i+1}=0$ to the product $w^+(k,R,T)$, while in the
second case, the $i$ contributes $(a_\lambda)_i$ for some
partition $\lambda$ occuring in the sequence $R$. Since there are
no occurrences of $i$ in $T$ between this $i$ and the relevant
marked $i$ in reverse column order, there cannot be any
occurrences of $\o{i}$ in $R$ between $\lambda$ and $\r{k-m-1}$,
so $(a_{\lambda})_i=(a_{\r{k-m-1}})_{i}=0$. Hence this tableau has
$w^+(k,R,T)$ and can safely be ignored. For the final two cases, a
similar argument shows that they give tableaux for which
$w^-(k,R,T)=0$.

Tableaux containing the following either cannot occur as they
cannot arise from a sequence of partitions or contribute zero to
the sum as in the above, so can also be ignored.

\[
\setlength{\unitlength}{17pt}
\begin{picture}(1,2)
\put(0,0){\line(1,0){1}} \put(0,0){\line(0,1){2}}
\put(0,1){\line(1,0){1}} \put(0,2){\line(1,0){1}}
\put(1,0){\line(0,1){2}} \put(0.4,1.3){$\overline{i}$}
\put(0.05,0.3){$\underline{\overline{i\!\!+\!\!1}}$}
\end{picture}
\setlength{\unitlength}{1pt}\qquad \setlength{\unitlength}{17pt}
\begin{picture}(1,2)
\put(0,0){\line(1,0){1}} \put(0,0){\line(0,1){2}}
\put(0,1){\line(1,0){1}} \put(0,2){\line(1,0){1}}
\put(1,0){\line(0,1){2}} \put(0.4,1.3){$\underline{\overline{i}}$}
\put(0.05,0.3){$\overline{i\!\!+\!\!1}$}
\end{picture}
\setlength{\unitlength}{1pt}\qquad \setlength{\unitlength}{17pt}
\begin{picture}(2,1) \put(0,0.5){\line(0,1){1}}
\put(0,0.5){\line(1,0){2}} \put(1,0.5){\line(0,1){1}}
\put(2,0.5){\line(0,1){1}} \put(0,1.5){\line(1,0){2}}
\put(0.4,0.75){$\underline{\overline{i}}$}
\put(1.4,0.75){$\overline{i}$}
\end{picture}
\setlength{\unitlength}{1pt} \qquad \setlength{\unitlength}{17pt}
\begin{picture}(2,1) \put(0,0.5){\line(0,1){1}}
\put(0,0.5){\line(1,0){2}} \put(1,0.5){\line(0,1){1}}
\put(2,0.5){\line(0,1){1}} \put(0,1.5){\line(1,0){2}}
\put(0.07,0.75){$\overline{i\!\!+\!\!1}$}
\put(1.07,0.75){$\underline{\overline{i\!\!+\!\!1}}$}
\end{picture}
\setlength{\unitlength}{1pt}
\]

This is because $\underline{\overline{i}}$ above $\overline{i+1}$
or $\overline{i}$ above $\overline{\underline{i+1}}$ cannot arise
from a sequence of partitions $R$. If $\underline{\overline{i}}$
to the immediate left of $\overline{i}$, then by examining $R$,
there would need to be a marked $i+1$ between these two entries in
reverse column order, which can only lie directly below
$\overline{i}$. But then to be semistandard, the entry directly
below $\underline{\overline{i}}$ must be an $i+1$. We have already
shown that this entry cannot be marked. So it is unmarked, in
which case we have already shown that the tableau gives zero
contribution so can be ignored. The case of
$\underline{\overline{i+1}}$ to the immediate right of
$\overline{i+1}$ proceeds similarly.

Wherever possible, we match up our tableaux as such:
\[
\begin{tabular}{|c|}
  \hline
  $\underline{\overline{i}}$ \\
  \hline
\end{tabular}+
\begin{tabular}{|c|}
  \hline
  $\underline{\overline{i}},i+1$ \\
  \hline
\end{tabular}\leftrightarrow
\begin{tabular}{|c|}
  \hline
  $\underline{\overline{i+1}}$ \\
  \hline
\end{tabular}+
\begin{tabular}{|c|}
  \hline
  $i,\underline{\overline{i+1}}$ \\
  \hline
\end{tabular}
\]
where all other elements of $T$ are unchanged. Then making this
pairing, the corresponding terms in $S_m$ cancel as it is a
special case of (\ref{theend}). Taking into account those tableaux
which we have already shown to give zero contribution to the sum,
we find that the only remaining tableaux $T$ for which the above
identification of positive and negative terms cannot be made, is
where $i$ and $i+1$, both barred and one underlined occur in the
same cell of $T$. We shall focus our attention on the positive
terms for which this happens, as the negative case proceeds
similarly.

Then the following must occur as a subsequence of $R$:

$$\r{k-m-2}\stackrel{i+1}\longrightarrow\r{k-m-1}\stackrel{i}\longrightarrow\r{k-m}\tangle\r{k}.$$

By replacing this with
$\r{k-m-2}\tangle\rho'\stackrel{i+1}\longrightarrow\r{k}$, we form
another sequence $R'$.

We map $T$ with $i$ and $i+1$ barred in the same box to $T'$ and
$T''$ with the barred $i$ removed in the first case and unbarred
in the second case.

\[
\begin{tabular}{|c|}
  \hline
  $\overline{i},\underline{\overline{i+1}}$ \\
  \hline
\end{tabular}\mapsto
\begin{tabular}{|c|}
  \hline
  $\underline{\overline{i+1}}$ \\
  \hline
\end{tabular}+\begin{tabular}{|c|}
  \hline
  $i,\underline{\overline{i+1}}$ \\
  \hline
\end{tabular}
\]

Again we have Lemma \ref{claim1},
\[
w^+(k,R,T)\Pi(\r{k-m})=\beta(w^-(k-1,R',T')+w^-(k-1,R',T''))\Pi(\r{k-m-2}).
\]
with the same proof.

Now we look at when the inverse to this map can be constructed:
\[
\begin{tabular}{|c|}
  \hline
  $\underline{\overline{i+1}}$ \\
  \hline
\end{tabular}+\begin{tabular}{|c|}
  \hline
  $i,\underline{\overline{i+1}}$ \\
  \hline
\end{tabular}\mapsto\begin{tabular}{|c|}
  \hline
  $\overline{i},\underline{\overline{i+1}}$ \\
  \hline
\end{tabular}
\]
Here we will underline the relevant $\overline{i+1}$ for clarity.

We consider the cases when this map cannot be made. They are

(i) Unbarred $i+1$ to the immediate left of
$\underline{\overline{i+1}}$: The unbarred $i+1$ contributes
$(a_{\r{k}})_{i+1}=0$ to the product, so these tableaux can be
ignored.

(ii). Unbarred $i$ above $\overline{\underline{i+1}}$: The
unbarred $i$ contributes  $(a_{\r{k-m-2}})_{i}=0$ to the product,
so again we have zero contribution so these tableaux can be
ignored. ignored.

(iii). $\overline{i+1}$ to the immediate left of
$\overline{\underline{i+1}}$: In order to have a sequence of
partitions, there must be a barred $i$ between these two barred
$i+1$'s, so must lie directly above the non-underlined one. Thus,
above $\overline{\underline{i+1}}$ must either lie an unbarred
$i$, giving case (ii), or a barred $i$, giving the final case
which we now deal with.

(iv). $\overline{i}$ immediately above
$\overline{\underline{i+1}}$: We note that in this case the
sequence must contain the subsequence
$\r{k-m-3}\stackrel{i}\longrightarrow\r{k-m-2}\tangle\rho'\stackrel{i+1}\longrightarrow\r{k}.$

Thus positive terms covered by Case 3 in $S_m$ give $\beta$ times
the negative terms in $\varepsilon_{m+1}$, with the exception of
those terms which are covered by case (iv) above.

Similarly, we see that the negative terms covered by Case 3 in
$S_m$ give $\beta$ times the positive terms in $\varepsilon_{m+1}$
with the same exception of those terms covered in case (iv).

Using our observation of the structure of case (iv) terms, we note
that if $(k,R,T)$ is such a positive $\varepsilon$-triple, then it
is also a negative $\varepsilon$-triple. Furthermore for such
terms $w^+(k,R,T)=w^-(k,R,T)$, so it is seen that such terms
cancel themselves out, giving no net contribution, and hence the
contribution to $S_m$ by all triples covered by this case is
exactly $-\beta\varepsilon_{m+1}$.

{\bf Case 4:} The only remaining triples are
$\varepsilon$-triples, so their contribution to $S_m$ is exactly
$\varepsilon_m$.

Hence we have proven (\ref{zem}), as desired. \end{proof}

\begin{remark} It is necessary for us to set some variables equal to zero
in (\ref{zsoln}), as otherwise, the formula does not hold true.
For example, if $n=2$, $\theta=\nu=(1^2)$ and $\mu=\phi$, then
using the recurrence, we calculate
$g_{\theta\mu}^\nu=\frac{1+\beta b_1}{1+\beta a_1}$ whereas
$\sum_{(R,T)}\beta^{|T|-|\theta|}w(T)=\frac{(1+\beta
b_1)^2}{(1+\beta a_1)(1+\beta a_2)}$. \end{remark}

\subsection{Specialisation to ordinary Grothendieck polynomials}

Specialisation to $a=0$ in (\ref{coef}) gives
$g_{\theta\mu}^\nu(0,0)=c_{\theta\mu}^\nu$, where $c_{\th\mu}^\nu$
is as defined in (\ref{deflr}). Under this specialisation, a pair
$(R,T)$ will contribute $\beta^{|T|-|\theta|}$ to the sum in
(\ref{zsoln}) if it consists entirely of barred entries, and 0
otherwise. For such tableaux, we must thus have $|T|=|\nu|-|\mu|$.
Hence we have the following.

\begin{theorem} \label{main} $c_{\theta\mu}^\nu$ is equal to
$\beta^{|\nu|-|\mu|-|\theta|}$ times the number of semistandard
set-valued $\theta$-tableaux with entries $r_1\prec
r_2\prec\ldots\prec r_l$ for which there is a related sequence of
partitions
$$\mu=\r{0}\o{1}\r{1}\o{2}\cdots\o{l}\r{l}=\nu.$$
\end{theorem}

We can show directly that this specialises to Buch's results as
quoted in Theorems \ref{b1} and \ref{b2}.

In the case $\theta=\lambda$, a partition, in order to directly
see that our result is equivalent to Theorem \ref{b1}, we note
that if $T$ is a $\mu$-tableau and $c(T)$ is to be a lattice word,
then the entries of $T$ are fixed, namely in that the $i$-th row
must contain only $i$'s. Call this particular tableau $T_\mu$. Now
for $c(T')$ to be a lattice word for some $\lambda*\mu$-tableau
$T'$, we must have $T'=T*T_\mu$ for some $\lambda$-tableau $T$.
Now we have a simple bijection between the two formulations in
this case, namely that which is given by $T\mapsto T*T_\mu$. The
condition of $c(T*T_\mu)$ being a lattice word is equivalent to
the sequence
$$\mu=\r{0}\o{1}\r{1}\o{2}\cdots\o{l}\r{l}=\nu.$$
consisting entirely of partitions.

For the case $\mu=\phi$, where we are expanding a skew
Grothendieck polynomial in the basis of ordinary Grothendieck
polynomials, Theorem \ref{b2} is easily seen to be consistent with
our formulation since $r_1,r_2,\ldots,r_m$ is a lattice word if
and only if
$$\phi=\r{0}\o{1}\r{1}\o{2}\cdots\o{m}\r{m}=\lambda$$ is a
sequence of partitions where $\lambda$ is the content of
$r_1,r_2,\ldots,r_m$.

\section{Grothendieck Polynomials via Isobaric Divided Differences}

The remainder of this paper will have a distinctly different flavour
to it, as we move away from calculating the Littlewood-Richardson
coefficients and instead devote the remainder of our energies to
exhibiting a relationship between the factorial Grothendieck
polynomials studied here, and the double Grothendieck polynomials,
as studied elsewhere. For the most part of this section, we follow
the exposition of Fomin and Kirillov \cite{FK}, supplying some proofs which are missing in their extended abstract.

\subsection{The symmetric group}

It is well known that the symmetric group $S_{n+1}$ is generated
by the $n$ simple reflections $s_i=(i,i+1)$, $i=1,2,\ldots,n$
subject to the relations
\begin{eqnarray*}
s_is_j&=&s_js_i\qquad\qquad\text{if\ } |i-j|\geq2,\\
s_is_{i+1}s_i&=&s_{i+1}s_is_{i+1},\\
s_i^2&=&1.
\end{eqnarray*}
For an element $w\in S_{n+1}$, let $\ell(w)$ denote the minimal
number $l$ for which $w$ can be written as a product of $l$ simple
reflections $w=s_{i_1}s_{i_2}\ldots s_{i_l}$. Then $\ell(w)=\#\{
i<j\mid w(i)>w(j) \}$. The longest word in $S_{n+1}$ is
$w_0=(n+1,n,\ldots,1)$ and satisfies $\ell(w_0)=\frac{n(n+1)}{2}$.

\subsection{Isobaric divided difference operators}

Let $R$ be a commutative ring with identity. We shall also assume
that $R$ contains various indeterminates used later, namely
$\beta$, $a_1,a_2,\ldots$. Let $f$ be a polynomial in the
variables $x_1,x_2,\ldots,x_{n+1}$ over $R$. For $i=1,2,\ldots,n$,
define the {\it isobaric divided difference operator} $\pi_i$ by
\[
\pi_if=\frac{(1+\beta
x_{i+1})f(\ldots,x_i,x_{i+1},\ldots)-(1+\beta
x_i)f(\ldots,x_{i+1},x_i,\ldots)}{x_{i}-x_{i+1}}
\]

Then these isobaric divided difference operators are easily
verified to satisfy the following relations:
\begin{eqnarray} \label{coxeter}
\pi_i\pi_j&=&\pi_j\pi_i\qquad\qquad\text{if\ } |i-j|\geq2,\\
\label{cox2}
\pi_i\pi_{i+1}\pi_i&=&\pi_{i+1}\pi_i\pi_{i+1},\\
\pi_i^2&=&-\beta \pi_i.
\end{eqnarray}

For each permutation $w\in S_{n+1}$ we now define the Grothendieck
polynomial $\GG_w$ in $x_1,x_2\ldots,x_{n+1},y_1,y_2,\ldots
y_{n+1}$. (We shall see later that these are actually polynomials
only in the variables $x_1,\ldots,x_n,y_1,\ldots,y_n$ but this
will not be immediately apparent.) If $w=w_0$, the longest
permutation, then set
$$\GG_{w_0}=\prod_{i+j\leq {n+1}}(x_i\oplus y_j).$$
If $w\neq w_0$, then there exists a simple reflection $s_i$ such
that $\ell(ws_i)>\ell(w)$. In such a case, we set
$\GG_w=\pi_i\GG_{ws_i}$. This definition is independent of the
choice of simple reflection, since the operators $\pi_i$ satisfy
the Coxeter relations (\ref{coxeter}) and (\ref{cox2}).

\subsection{The algebra $H_n$.}

Let $H_n$ be the $R[x_1,\ldots,x_n]$-algebra generated by
$u_1,u_2,\ldots ,u_n$ subject to
\begin{eqnarray}
\label{defhn1} u_iu_j&=&u_ju_i\qquad\qquad\text{if\ } |i-j|\geq2,\\
u_iu_{i+1}u_i&=&u_{i+1}u_iu_{i+1},\\
\label{defhn3} u_i^2&=&\beta u_i.
\end{eqnarray}

Then $H_n$ has dimension $(n+1)!$ and a natural basis $u_w$
indexed by elements of $S_{n+1}$ where $u_w=u_{i_1}\ldots u_{i_l}$
if $w=s_{i_1}\ldots s_{i_l}$ is a minimal representation of $w$ as
a product of simple reflections. \cite[Ch 4, \S2, Ex 23]{bourbaki}

For $x\in R[x_1,\ldots,x_n]$, define the following:
\begin{eqnarray*}
h_i(x)&=&1+xu_i \\
 A_i(x)&=&h_{n}(x)h_{n-1}(x)\ldots
h_i(x), \\
 B_i(x)&=&h_i(x)h_{i+1}(x)\ldots h_{n}(x), \\ A(x)&=&A_1(x), \\
 B(x)&=&B_1(x).
 \end{eqnarray*}

 For $x=\vec{x}$, we also define
\begin{eqnarray*}
\GG(x)&=&A_1(x_1)A_2(x_2)\ldots A_{n}(x_{n}),
\\\overline\GG(x)&=&B_{n}(x_{n})B_{n-1}(x_{n-1})\ldots
B_1(x_1).\end{eqnarray*}

Now we begin proving some preliminary identities in $H_n$.

\begin{lemma} \label{lemma}\cite{FK}
\begin{equation*}
B(x)B(y)=B(y)B(x).
\end{equation*}\end{lemma}
\begin{proof}
Expand $B(x)B(y)$ as a sum of $4^n$ terms. We identify each of
these terms with a 2-colouring of a $2\times n$ array of boxes.
The two colours chosen here are crossed and uncoloured. Each
colouring $U$ is identified with the term obtained by taking
$xu_i$ from the factor $h_i(x)$ (respectively $yu_i$ from
$h_i(y)$) from the $i$-th box in the first (respectively second)
row if it contains a cross and $1$ otherwise. Denote this term by
$\xi(U)$. So for example if
\[
\setlength{\unitlength}{10pt} U=\
\parbox[c]{5.2\unitlength}{\begin{picture}(5,2)
\multiput(0,0)(0,1){3}{\line(1,0){5}}
\multiput(0,0)(1,0){6}{\line(0,1){2}}
\multiput(0,1)(2,0){3}{\line(1,1){1}}
\multiput(0,2)(2,0){3}{\line(1,-1){1}}
\multiput(2,0)(1,0){2}{\line(1,1){1}}
\multiput(2,1)(1,0){2}{\line(1,-1){1}}
\end{picture}} \setlength{\unitlength}{1pt},
\]
then $\xi(U)=x^3y^2u_1u_3u_5u_3u_4$.

Let $\sigma(x^i y^j u_w)=x^j y^i u_w$.

We now find a bijection $U\to U'$ such that
$\sigma(\xi(U))=\xi(U')$, assuming inductively, that such a
bijection exists for all such $2\times m$ arrays with $m<n$.

There are four different types of columns that can occur in $U$,
which we shall unimaginatively call types I, II, III and IV as
follows:
\[
\rm I:\ \setlength{\unitlength}{10pt}
\begin{picture}(2,2)
\put(0,0){\line(1,0){1}} \put(0,0){\line(0,1){2}}
\put(0,1){\line(1,0){1}} \put(0,2){\line(1,0){1}}
\put(1,0){\line(0,1){2}}
\end{picture}
\setlength{\unitlength}{1pt} II:\ \setlength{\unitlength}{10pt}
\begin{picture}(2,2)
\put(0,0){\line(1,0){1}} \put(0,0){\line(0,1){2}}
\put(0,1){\line(1,0){1}} \put(0,2){\line(1,0){1}}
\put(1,0){\line(0,1){2}} \put(0,1){\line(1,1){1}}
\put(1,1){\line(-1,1){1}}
\end{picture}
\setlength{\unitlength}{1pt} III:\ \setlength{\unitlength}{10pt}
\begin{picture}(2,2)
\put(0,0){\line(1,0){1}} \put(0,0){\line(0,1){2}}
\put(0,1){\line(1,0){1}} \put(0,2){\line(1,0){1}}
\put(1,0){\line(0,1){2}} \put(0,0){\line(1,1){1}}
\put(1,0){\line(-1,1){1}}
\end{picture}
\setlength{\unitlength}{1pt} IV:\ \setlength{\unitlength}{10pt}
\begin{picture}(2,2)
\put(0,0){\line(1,0){1}} \put(0,0){\line(0,1){2}}
\put(0,1){\line(1,0){1}} \put(0,2){\line(1,0){1}}
\put(1,0){\line(0,1){2}} \put(0,0){\line(1,1){1}}
\put(1,0){\line(-1,1){1}} \put(0,1){\line(1,1){1}}
\put(1,1){\line(-1,1){1}}
\end{picture}
\setlength{\unitlength}{1pt}
\]

{\bf Case:} $U$ contains a type I column.

\begin{eqnarray*}
\hbox{Suppose}\hspace{13pt} \setlength{\unitlength}{10pt} U&=&\
\parbox[c]{6\unitlength}{\begin{picture}(7,2)
\multiput(0,0)(0,2){2}{\line(1,0){7}} \put(3,1){\line(1,0){1}}
\put(0,0){\line(0,1){2}}\put(3,0){\line(0,1){2}}\put(4,0){\line(0,1){2}}\put(7,0){\line(0,1){2}}
\put(1,0.7){$U_1$} \put(5,0.7){$U_2$}
\end{picture}} \setlength{\unitlength}{1pt} \\
\hbox{Define }\hspace{16pt} \setlength{\unitlength}{10pt} U'&=&\
\parbox[c]{6\unitlength}{\begin{picture}(7,2)
\multiput(0,0)(0,2){2}{\line(1,0){7}} \put(3,1){\line(1,0){1}}
\put(0,0){\line(0,1){2}}\put(3,0){\line(0,1){2}}\put(4,0){\line(0,1){2}}\put(7,0){\line(0,1){2}}
\put(1,0.7){$U'_1$} \put(5,0.7){$U'_2$}
\end{picture}} \setlength{\unitlength}{1pt}
\end{eqnarray*}

Terms from the lower half of $U_1$ and upper half of $U_2$ always
commute. So as $\sigma(\xi(U_i))=\xi(U_i')$ for $i=1,2$,
$\sigma(\xi(U))=\xi(U)$ as required.

So now we may assume that $U$ contains no type I columns.

{\bf Case:} $U$ contains a type III column to the left of a type
II column.

Of all such pairs of columns, consider a minimally separated pair.
Then $U$ must contain only type IV columns between this minimally
separated pair, so must be of the form
\begin{eqnarray*}
 \setlength{\unitlength}{10pt} U&=&\
\parbox[c]{6\unitlength}{\begin{picture}(11,2)
\multiput(0,0)(0,2){2}{\line(1,0){11}}
\put(3,1){\line(1,0){1}}\put(7,1){\line(1,0){1}}
\put(0,0){\line(0,1){2}}\put(3,0){\line(0,1){2}}\put(4,0){\line(0,1){2}}\put(7,0){\line(0,1){2}}
\put(8,0){\line(0,1){2}}\put(11,0){\line(0,1){2}}
\put(1,0.7){$U_1$} \put(9,0.7){$U_2$}
\put(3,0){\line(1,1){1}}\put(4,0){\line(1,1){1}}\put(6,0){\line(1,1){1}}
\put(4,1){\line(1,1){1}}\put(6,1){\line(1,1){1}}\put(7,1){\line(1,1){1}}
\put(3,1){\line(1,-1){1}} \put(4,2){\line(1,-1){1}}
\put(6,1){\line(1,-1){1}} \put(6,2){\line(1,-1){1}}
\put(4,1){\line(1,-1){1}} \put(7,2){\line(1,-1){1}}
\put(4.99,0.98){$\ldots$}
\end{picture}} \setlength{\unitlength}{1pt}\qquad\qquad\qquad  \\
\hbox{Define }\hspace{16pt} \setlength{\unitlength}{10pt} U'&=&\
\parbox[c]{6\unitlength}{\begin{picture}(11,2)
\multiput(0,0)(0,2){2}{\line(1,0){11}}
\put(3,1){\line(1,0){1}}\put(7,1){\line(1,0){1}}
\put(0,0){\line(0,1){2}}\put(3,0){\line(0,1){2}}\put(4,0){\line(0,1){2}}\put(7,0){\line(0,1){2}}
\put(8,0){\line(0,1){2}}\put(11,0){\line(0,1){2}}
\put(1,0.7){$U_1'$} \put(9,0.7){$U_2'$}
\put(3,0){\line(1,1){1}}\put(4,0){\line(1,1){1}}\put(6,0){\line(1,1){1}}
\put(4,1){\line(1,1){1}}\put(6,1){\line(1,1){1}}\put(7,1){\line(1,1){1}}
\put(3,1){\line(1,-1){1}} \put(4,2){\line(1,-1){1}}
\put(6,1){\line(1,-1){1}} \put(6,2){\line(1,-1){1}}
\put(4,1){\line(1,-1){1}} \put(7,2){\line(1,-1){1}}
\put(4.99,0.98){$\ldots$}
\end{picture}} \setlength{\unitlength}{1pt}\
\end{eqnarray*}

Terms from the lower half of $U_1$ commute with those from the
upper part of $U$ to the right of $U_1$, and terms from the upper
half of $U_2$ commute with those from the lower half of $U$ to the
left of $U_2$. So, again using strong induction, from
$\sigma_i(\xi(U_i))=\xi(U_i')$, we get $\sigma(\xi(U))=\xi(U_i')$
as required.

So now may also assume that such an arrangement does not exist.
Thus, we are only left to consider $U$ of the following schematic
type:

\[
\setlength{\unitlength}{13pt}
\begin{picture}(10,2)
\multiput(0,2)(0,-1){3}{\line(1,0){10}}  \put(0,0){\line(0,1){2}}
\put(1,0){\line(0,1){2}} \put(2,0){\line(0,1){2}}
\put(0,1){\line(1,1){1}} \put(1,1){\line(-1,1){1}}
\put(1,1){\line(1,1){1}} \put(2,1){\line(-1,1){1}}
\put(2.5,0.45){$\ldots$} \put(2.5,1.45){$\ldots$}
\multiput(4,0)(1,0){4}{\line(0,1){2}} \put(5,0){\line(0,1){2}}
\put(4,1){\line(1,1){1}} \put(5,1){\line(-1,1){1}}
\put(5,0){\line(1,1){1}} \multiput(6,0)(1,0){2}{\line(-1,1){1}}
\put(7.5,0.45){$\ldots$} \put(7.5,1.45){$\ldots$}
\put(10,0){\line(0,1){2}} \put(6,0){\line(1,1){1}}
\put(9,0){\line(0,1){2}} \put(9,0){\line(1,1){1}}
\put(10,0){\line(-1,1){1}} \put(0.35,0.2){$?$} \put(1.35,0.2){$?$}
\put(4.35,0.2){$?$} \put(9.35,1.2){$?$} \put(6.35,1.2){$?$}
\end{picture}
\setlength{\unitlength}{1pt}
\]

Suppose that $U$ contains $a$ type II columns and $b$ type III
columns. Without loss of generality, let us assume that $a\geq b$.
We can do so, since the $a<b$ case proceeds similarly, or
alternatively and equivalently can be defined as the inverse of
the $a>b$ case.

Let $V$ denote that part of $U$ lying strictly between the $b$-th
type II column (counting from the left) and the leftmost type III
column.

Draw a horizontal cutting line through the middle of $V$. Now draw
a vertical cutting line one boxwidth from the right hand edge of
$V$. Glue together opposite edges of $V$ to form a torus, and cut
this torus along the cutting lines constructed, while preserving
the directional notions of up, down, left and right to create
$V'$.

So for example if
\[
\setlength{\unitlength}{10pt} V=\
\parbox[c]{5.2\unitlength}{\begin{picture}(5,2)
\multiput(0,0)(0,1){3}{\line(1,0){5}}
\multiput(0,0)(1,0){6}{\line(0,1){2}}
\multiput(0,1)(1,0){5}{\line(1,1){1}}
\multiput(0,2)(1,0){5}{\line(1,-1){1}}
\multiput(3,0)(1,0){2}{\line(1,1){1}}
\multiput(3,1)(1,0){2}{\line(1,-1){1}} \put(1,0){\line(1,1){1}}
\put(1,1){\line(1,-1){1}}
\end{picture}} \setlength{\unitlength}{1pt},
\text{ then  \ } \setlength{\unitlength}{10pt} V'=\
\parbox[c]{5.2\unitlength}{\begin{picture}(5,2)
\multiput(0,0)(0,1){3}{\line(1,0){5}}
\multiput(0,0)(1,0){6}{\line(0,1){2}}
\multiput(0,1)(2,0){3}{\line(1,1){1}}
\multiput(0,2)(2,0){3}{\line(1,-1){1}}
\multiput(0,0)(1,0){5}{\line(1,1){1}}
\multiput(0,1)(1,0){5}{\line(1,-1){1}}
\end{picture}} \setlength{\unitlength}{1pt}.
\]

Returning to $U$, construct $U'$ by replacing $V$ with $V'$. It is
clear that the degrees of $x$ and $y$ in $\sigma(\xi(U))$ and
$\xi(U')$ match, while the fact that the same basis element $u_w$
is obtained in each case follows from the following relations:
\begin{eqnarray*}
(u_i u_{i+1}\ldots u_k)u_j&=&u_{j+1}(u_i u_{i+1}\ldots u_k)\qquad
(i\leq j<k)\\
(u_i u_{i+1}\ldots u_k)u_k&=&u_i(u_i u_{i+1}\ldots u_k)
\end{eqnarray*}

Thus we have our bijection $U\mapsto U'$ such that
$\sigma(\xi(U))=\xi(U')$ as required. Hence,

\[
B(x)B(y)=\sum_U \xi(U)=\sum_{U'}\xi(U')=\sum_U
\sigma(\xi(U))=B(y)B(x)
\]
and the lemma is proven.

\end{proof}

\begin{corollary} \label{corollary}\cite{FK}
$A_i(x)$ and $A_i(y)$ commute, as do $A(x)$ and $A(y)$.
\end{corollary}

The following lemma is easily proven by expanding out and applying
the defining relations (\ref{defhn1}) to (\ref{defhn3}). These
three simple equations are used extensively in the following work.

\begin{lemma}\cite{FK}
\begin{eqnarray}
h_i(x)h_j(y)&=&h_j(y)h_i(x)\qquad\text{if\ \ }|i-j|\geq 2\\
h_i(x)h_{i+1}(x\oplus y)h_i(y)&=&h_{i+1}(y)h_i(x\oplus y)h_{i+1}(x)\\
h_i(x)h_i(y)&=&h_i(x\oplus y)
\end{eqnarray}
\end{lemma}

\begin{lemma} \label{big} \cite{fkmini,FK}
\[
(\pi_i+\beta)\GG(x)=\GG(x) u_i
\]
\end{lemma}
(Here, $\pi_i$ is acting only on the elements of our coefficient
ring $R[x_1,\ldots,x_n]$, so we could say that $\pi_i$ acts
trivially on $u_j$ for all $j$.)

\begin{proof}
Write
\begin{equation*}
\GG(x)=A_1(x_1)\ldots A_i(x_i)A_i(x_{i+1})h_{i}(\ominus
x_{i+1})A_{i+2}(x_{i+2})\ldots A_{n}(x_{n}).
\end{equation*}
A routine calculation shows that $(\pi_i+\beta)h_{i}(\ominus
x_{i+1})=h_i(\ominus x_{i+1})u_i$. By Corollary \ref{corollary},
$A_i(x_i)$ and $A_i(x_{i+1})$ commute so $A_1(x_1)\ldots
A_i(x_i)A_{i}(x_{i+1})$ is symmetric in $x_i$ and $x_{i+1}$. Since
$\pi_i(fg)=f\pi_i g$ whenever $f$ is symmetric in $x_i$ and
$x_{i+1}$, we have
\begin{eqnarray*}
&&\hspace{-25pt}(\pi_i+\beta)\GG(x) \\&=&A_1(x_1)\ldots
A_i(x_i)A_{i}(x_{i+1})(\pi_i+\beta)h_{i}(\ominus
x_{i+1})A_{i+2}(x_{i+2})\ldots A_n(x_n) \\
&=&A_1(x_1)\ldots A_i(x_i)A_{i}(x_{i+1})h_{i}(\ominus
x_{i+1})u_i A_{i+2}(x_{i+2})\ldots A_n(x_n) \\
&=&\GG(x)u_i.
\end{eqnarray*} as required since $u_i$ and $A_j(x_j)$ commute for
$j\geq i+2$.
\end{proof}

\begin{lemma}\cite{FK}
$A_i(x)$ and $B_i(y)$ commute.
\end{lemma}

\begin{proof} We prove this by descending induction on $i$.
\begin{eqnarray*}
A_i(x)B_i(y)&=&h_n(x)\ldots h_{i+1}(x)h_i(x\oplus
y)h_{i+1}(y)\ldots
h_n(y) \\
&=& h_n(x)\ldots h_{i+2}(x)h_i(y)h_{i+1}(x\oplus
y)h_i(x)h_{i+2}(y)\ldots h_n(y) \\
&=& h_i(y)h_n(x)\ldots
h_{i+2}(x)h_{i+1}(x)h_{i+1}(y)h_{i+2}(y)\ldots h_n(y) \\
&=& h_i(y) A_{i+1}(x)B_{i+1}(y) h_i(x) \\
&=& h_i(y) B_{i+1}(y)A_{i+1}(x)h_i(x) \\
&=& B_i(y)A_i(x) \qquad \ \hbox{as required.}
\end{eqnarray*}
\end{proof}

\begin{lemma}\label{5.2} \cite{FK}
\[
B_{n}(y_{n})\ldots B_i(y_i)A_i(x)=h_{n}(x\oplus y_{n})\ldots
h_i(x\oplus y_i)B_{n}(y_{n-1})\ldots B_{i+1}(y_i).
\] \end{lemma}
\begin{proof}
Again, we use descending induction on $i$.
\begin{eqnarray*}
\hbox{LHS}&=&B_n(y_n)\ldots B_{i+1}(y_{i+1})A_i(x)B_i(y_i) \\
&=& B_n(y_n)\ldots B_{i+1}(y_{i+1})A_{i+1}(x)h_i(x\oplus
y_i)B_{i+1}(y_i) \\
&=& h_n(x\oplus y_{n})\ldots h_{i+1}(x\oplus
y_{i+1})B_{n}(y_{n-1})\ldots B_{i+2}(y_{i+1})h_i(x\oplus
y_i)B_{i+1}(y_i) \\
&=& \hbox{RHS}
\end{eqnarray*}
since $h_{i}(x\oplus y_i)$ commutes with $B_j(y_{j-1})$ for $j\geq
i+2$.
\end{proof}

\begin{lemma}\label{l5.3}\cite{FK}
\begin{equation}\label{5.3}
\overline\GG(y)\GG(x)=\prod_{i=1}^n \prod_{j=n+1-i}^1
h_{i+j-1}(x_i\oplus y_j).
\end{equation}
\end{lemma}
\begin{proof}
We prove this result by induction on $n$. Lemma \ref{5.2} gives
$$\overline\GG(y)\GG(x)=h_n(x_1\oplus y_n)\ldots h_1(x_1\oplus
y_1)B_n(y_{n-1})\ldots B_2(y_1)A_2(x_2)\ldots A_n(x_n)$$ which
equals our desired result by applying the inductive hypothesis.
\end{proof}

\subsection{A generating function for Grothendieck polynomials}

\begin{theorem}\cite{FK}
\begin{equation} \label{9.12} \overline\GG(y)\GG(x)=\sum_{w\in
S_{n+1}}\GG_wu_w.
\end{equation}
\end{theorem}

\begin{proof}
Let $g_w$ be the coefficient of $u_w$ in $\overline\GG(y)\GG(x)$.
We prove by decreasing induction on $\ell(w)$ that $g_w=\GG_w$.

First we consider the case $w=w_0$. There are
$\frac{n(n+1)}{2}=\ell(w_0)$ terms in the product on the right
hand side of (\ref{5.3}). We also have $$(u_n\ldots u_1)(u_n\ldots
u_2)\ldots(u_nu_{n-1})(u_n)=u_{w_0}.$$ Hence,
$$g_{w_0}=\prod_{i=1}^n\prod_{j=n-i+1}^1 x_i\oplus
y_j=\prod_{i+j\leq n+1} x_i\oplus y_j =\GG_{w_0}.$$

Now suppose $w\neq w_0$ and consider a simple reflection $s_i$
such that $\ell(ws_i)>\ell(w_i)$. Lemma \ref{big} tells us that
$(\pi_i+\beta)\overline\GG(y)\GG(x)=\overline\GG(y)\GG(x)u_i$.
Comparing the coefficient of $u_{ws_i}$ in this equation then
gives $$\pi_ig_{ws_i}+\beta g_{ws_i}=g_w+\beta g_{ws_i}.$$ By our
inductive assumption, $g_{ws_i}=\GG_{ws_i}$, so $g_w=\pi_i
\GG_{ws_i}=\GG_w$ as required and the theorem is proven.
\end{proof}

Combining (\ref{5.3}) and (\ref{9.12}) gives
\[
\prod_{i=1}^{n+m}\prod_{j=n+m-i}^{1} h_{i+j-1}(x_i\oplus
y_j)=\sum_{w\in S_{n+m}} \GG_w u_w.
\]

If we now temporarily restrict ourselves to the finite set of
variables $x_1,x_2,\ldots,x_k,y_1,y_2,\ldots,y_l$ (by setting
$x_i=0$ if $i>k$ and $y_j=0$ if $j>l$), then for $m>\max(k,l)$, we
apply the homomorphism $\psi:H_{n+m}\to H_n$ given by
$\psi(h_i)=0$ if $i\leq m$ and $\psi(h_i)=h_{i-m}$ if $i>m$ to get
\[
B(y_l)B(y_{l-1})\ldots B(y_1)A(x_1)A(x_2)\ldots A(x_k)=\sum_{w\in
S_{n+1}}\GG_{1^m\times w}u_w.
\]
Here, for $w\in S_{n+1}$, $1^m\times w\in S_{m+n+1}$ is the
permutation with $(1^m\times w)(i)=i$ if $i\leq m$ and $(1^m\times
w)(i)=m+w(i-m)$ if $i>m$.

Thus the coefficient of each fixed monomial in $\GG_{1^m\times w}$
eventually becomes stable as $m$ tends to infinity. So now we can
make the following definition:

\begin{definition}
For a permutation $w$, the double stable Grothendieck polynomial
in $x$ and $y$ is defined to be the power series
$$G_w(x;y)=\lim_{m\to\infty}\GG_{1^m\times w}.$$
\end{definition}

Restricting ourselves again to the finite set of variables
$x_1,\ldots x_k,y_1\ldots y_l$, we thus have
\begin{equation} \label{9.n}
B(y_l)B(y_{l-1})\ldots B(y_1)A(x_1)A(x_2)\ldots A(x_k)=\sum_{w\in
S_{n+1}}G_{w}(x;y)u_w.
\end{equation}

\begin{lemma}
Let $p$ be an integer. Then
\begin{equation} \label{nearly}
\prod_{m=\infty}^{-\infty}\prod_{i=n}^1 h_i(x_m\oplus
y_{m+i-p})=\sum_{w\in S_{n+1}} G_w(x;y)u_w
\end{equation}
where any out of range variables are set equal to zero.
\end{lemma}

\begin{proof}
Repeated application of Lemma \ref{5.2} shows the left hand side
of (\ref{nearly}) to equal $$B(y_l)B(y_{l-1})\ldots
B(y_1)A(x_k)A(x_{k-1})\ldots A(x_1).$$ To complete the proof of
the lemma, we use Corollary \ref{corollary} which tells us that
$A(x_i)$ and $A(x_j)$ commute, together with (\ref{9.n}) and we
are done.
\end{proof}

\section{Relationship between factorial and double Grothen\-dieck
Polynomials}

We shall restrict ourselves now to considering factorial
Grothendieck polynomials $G_\lambda(x|a)$ for $\lambda$ a
partition. By making such a restriction, rather than considering
$G_\th(x|a)$ for an arbitrary skew diagram $\theta$, it enables
our main result in this section, namely Theorem \ref{final} to be
stated in simple terms. However, one can define a double
Grothendieck polynomial $G_\th(x;y)$ for a skew partition $\theta$
in a similar vein, as appears in \cite{B:buch}. A natural
conjecture would be that a result similar to Theorem \ref{final}
should exist relating $G_\th(x;y)$ and $G_\theta(x|y)$ for a skew
partition $\th$.

First, we need a preliminary definition before we can define the
double Grothendieck polynomial $G_\lambda(x;y)$ for a partition
$\lambda$. So suppose
$\lambda=(\lambda_1,\lambda_2,\ldots\lambda_p)$ is a partition.
Here, we do not necessarily have $p=\ell(\lambda)$, but certainly
we must have $p\geq\ell(\lambda)$. Define the permutation
$w(\lambda)\in S_\infty=\varinjlim{S_n}$ by
$w(\lambda)(i)=i+\lambda_{p+1-i}$ for $1\leq i\leq p$ and
$w(\lambda)(i)=i-\lambda'_{i-p}$ for $i>p$. An explicit
representation of this permutation as a product of simple
reflections is constructed in Lemma \ref{perm}, so $w(\lambda)$ is indeed a permutation. Now we can make our definition:

\begin{definition} Define the double Grothendieck polynomial
$G_\lambda(x;y):=G_{w(\lambda)}(x;y).$ Note that this definition
is independent of $p$, since $G_w(x;y)=G_{1^m\times w}(x;y)$.
\end{definition}

We now proceed in a similar vein to Buch \cite{B:buch}.

Place a {\it diagonal numbering} in the boxes of $\lambda$ as
follows: Number the NW-SE (defining compass directions north,
west, south and east on $\lambda$ in the usual manner so that
north is at the top and west is on the left) diagonals of
$\lambda$ with positive integers, consecutively increasing from SW
to NE, such that $p$ is the number of the main diagonal. For
example with $\lambda=(4,3,1)$ and $p=4$ the numbering is
explicitly shown in the following picture:

\[
\setlength{\unitlength}{1.4em}
\begin{picture}(4,3)
\multiput(0,1)(0,1){3}{\line(1,0){3}}
\multiput(0,1)(1,0){4}{\line(0,1){2}}
\multiput(0,0)(3,2){2}{\line(1,0){1}}
\multiput(0,0)(1,0){2}{\line(0,1){1}} \put(4,3){\line(-1,0){1}}
\put(4,3){\line(0,-1){1}} \put(0.35,2.25){4} \put(0.35,1.25){3}
\put(0.35,0.25){2} \put(1.35,2.25){5} \put(1.35,1.25){4}
\put(2.35,2.25){6} \put(2.35,1.25){5} \put(3.35,2.25){7}
\end{picture}
\setlength{\unitlength}{1pt}
\]


Say that a partition $\mu$ contains an outer corner in the $i$-th
diagonal if this diagonal contains a box outside $\mu$ such that
the two boxes immediately above and to the left of it are in
$\mu$. Say that $\mu$ contains an inner corner in the $i$-th
diagonal if this diagonal contains a box inside $\mu$ such that
the two boxes immediately below and to the right of it are not in
$\mu$. So continuing with the example above, $\lambda$ contains an
outer corner (among others) in the third diagonal and an inner
corner in the fifth diagonal.

Suppose that $n$ is such that $n\geq p+\lambda_1-1$ and let
$V=\bigoplus_{\mu}R[x_1,\ldots,x_n]\cdot[\mu]$ be the free
$R[x_1,\ldots,x_n]$-module with basis indexed by partitions $\mu$.
As in \cite{B:buch}, we define an action of $H_n$ on $V$ as
follows:

If $\mu$ has an outer corner in the $i$-th diagonal, set
$u_i[\mu]=[\tilde\mu]$ where $\tilde\mu$ is the partition obtained
from $\mu$ by adding a box in this corner. If $\mu$ has an inner
corner in the $i$-th diagonal, set $u_i[\mu]=\beta[\mu]$. In all
other cases, set $u_i[\mu]=0$.





\begin{lemma} \label{perm}
We have the following representation of $w(\lambda)$ as a product
of simple reflections. Let $(i_1,i_2,\ldots ,i_{|\lambda|})$ be
the diagonal numbers of the boxes of $\lambda$, read one row at a
time from bottom to top, reading from right to left in each row.
Then
$w(\lambda)=s_{i_1}s_{i_2}\ldots s_{i_{|\lambda|}}$. Furthermore,
this representation is minimal, that is
$\ell(w(\lambda))=|\lambda|$, we have the identity $u_{w(\lambda)}
[\phi]=[\lambda]$ and in any expression of the form
$\beta^{t-|\lambda|}u_{w(\lambda)}=u_{i_1} u_{i_2}\ldots u_{i_t}$,
we have $i_t=p$.
\end{lemma}

\begin{proof}

Suppose that $i\leq p$. Then in calculating $s_{i_1}s_{i_2}\ldots
s_{i_{|\lambda|}}(i)$, the relevant simple reflections are exactly
those in the $(p+1-i)$-th row of $\lambda$, so
$s_{i_1}s_{i_2}\ldots s_{i_{|\lambda|}}(i)=i+\lambda_{p+1-i}$. If
$i>p$, then the relevant simple reflections are exactly those in
the $(i-p)$-th column of $\lambda$, so $s_{i_1}s_{i_2}\ldots
s_{i_{|\lambda|}}(i)=i-\lambda'_{i-p}$ in this case. Hence
$w(\lambda)=s_{i_1}s_{i_2}\ldots s_{i_{|\lambda|}}$.

Note that $w(\lambda)(i+1)>w(\lambda)(i)$ for all $i\neq p$. Hence
if $i<j$ is such that $w(\lambda)(i)>w(\lambda)(j)$, we must have
$i\leq p$ and $j>p$. For $i$ and $j$ in this range,
$w(\lambda)(i)>w(\lambda)(j)$ if and only if
$1+\lambda_{p+1-i}+\lambda'_{j-p}>(p+1-i)+(j-p)$, which occurs if
and only if $(p+1-i,j-p)\in\lambda$. Hence there are $|\lambda|$
such pairs $(i,j)$, so $\ell(w(\lambda))=|\lambda|$.

For the remaining statement, we first notice that since
$w(\lambda)=s_{i_1}\ldots s_{i_{|\lambda|}}$, we easily calculate
$u_{w(\lambda)} [\phi]=[\lambda]$. Now if
$\beta^{t-|\lambda|}u_{w(\lambda)}=u_{i_1} u_{i_2}\ldots u_{i_t}$,
then $u_{i_1} u_{i_2}\ldots
u_{i_t}[\phi]=\beta^{t-|\lambda|}u_{w(\lambda)}[\phi]=\beta^{t-|\lambda|}[\lambda]\neq
0$. Hence $u_{i_t}[\phi]\neq 0$, so $i_t=p$.

\end{proof}

\begin{lemma}\cite{B:buch} If $w\in S_{n+1}$ is such that $u_w[\phi]\neq 0$,
then $u_w[\phi]=[\mu]$ for some partition $\mu$ and furthermore
$w=w(\mu)$.
\end{lemma}

\begin{proof}
The proof of this lemma by induction on $\ell(w)$, and contained
in \cite{B:buch}. \end{proof}




We have the following immediate corollary:

\begin{corollary}
For all $h\in H_n$, the coefficient of $u_{w(\lambda)}$ in $h$ is
equal to the coefficient of $[\lambda]$ in $h[\phi]$.
\end{corollary}

Define the products $P$ and $Q$ by
\[
P=\prod_{m=\infty}^{1}\prod_{i=n}^1 h_i(x_m\oplus y_{m+i-p}).
\]

\[
Q=\prod_{m=k}^1\prod_{i=n}^1 h_i(x_m\oplus y_{m+i-p}).
\]

\begin{theorem}
The coefficient of $u_{w(\lambda)}$ in $P$ is the double
Grothendieck polynomial $G_\lambda(x;y)$.
\end{theorem}

\begin{proof}
This follows from Lemma \ref{9.12} and Lemma \ref{perm}, noting
that factors on the left hand side of (\ref{nearly}) with $m\leq
0$ are either one, or do not contain $u_p$. \end{proof}

\begin{theorem}
The coefficient of $[\lambda]$ in $Q[\phi]$ is the factorial
Grothendieck polynomial $G_\lambda(x|y)$ in $x_1,\ldots x_k$.
\end{theorem}

\begin{proof}
Expand $Q$, and note that each term is a product of terms of the
form $(x_m\oplus y_{m+i-p} )u_i$. Write this product as
$$\prod_{j=1}^q(x_{m_j}\oplus y_{m_j+i_j-p})u_{i_j}.$$ If
$(\prod_{j=1}^q u_{i_j} )[\phi]\neq 0$, then we can interpret this
product in the following way:

Form the tableau $T$ by placing $m_j$ in the inner corner in
diagonal number $i_j$ of the partition $u_{i_j}\ldots
u_{i_q}[\phi]$ for $j=1,2,\ldots q$. These numbers are added in
nondecreasing order, and the occurrences of each number $i$ are
added from SW to NE. Furthermore, at all stages during the
addition process, the shape of all the numbers added up to that
point is a partition. So $T$ is a semistandard set-valued tableau
with entries from $[n]$. Note that if $\alpha$ is a cell with
diagonal number $i$, then $c(\alpha)=i-p$. Hence we can write
$$\Big(\prod_{j=1}^q u_{i_j} \Big)[\phi]= \beta^{|T|-|\lambda|}
\Big(\prod_{\substack{\alpha\in\lambda \\ r\in
T(\alpha)}}x_r\oplus y_{r+c(\alpha)}\Big)
 u_{w(\lambda)}$$ for
some partition $\lambda$ for which $T$ is a $\lambda$-tableau.

Upon considering all such terms in $Q$, it becomes evident that
the coefficient of $[\lambda]$ in $Q[\phi]$ is the factorial
Grothendieck polynomial $G_\lambda(x|y)$ as required.
\end{proof}

Note that $\lim_{k\to\infty}Q=P$. Hence, the preceding three
results give us the following theorem.

\begin{theorem}[Relationship between factorial and double
Grothendieck polynomials]\label{final}
\[
G_\lambda(x;y)=\lim_{k\to\infty}G_\lambda(x_1,\ldots,x_k|y).
\]
\end{theorem}

\section*{Acknowledgements}
The author would like to thank A. Molev for his guidance through this area of
study, and also A. Henderson, for his help with preparation of the manuscript.

\bibliographystyle{amsplain}

\begin{thebibliography}{10}

\bibitem[Bo]{bourbaki}
Bourbaki, N, {\em``Groupes et alg\`ebres de Lie,"} Ch 4-6, Hermann,
Paris 1968.

\bibitem[Bu]{B:buch}
{A. S. Buch}, {\em A Littlewood-Richardson Rule for the $K$-Theory
of Grassmannians},
   Acta. Math, {\bf 189} (2002), 37--78.

\bibitem[FG]{FG}
{S. Fomin and C. Greene}. {\em Noncommutative Schur functions and
their applications}, Discrete Math, {\bf 193}, (1998), 179--200.

\bibitem[FK1]{fkmini}
{S. Fomin and A. N. Kirillov}, {\em Grothendieck polynomials and
the Yang-Baxter equation}, Proc. Formal Power Series and Alg.
Comb, (1994), 183--190.

\bibitem[FK2]{FK}
{S. Fomin and A. N. Kirillov}, {\em Yang-Baxter equation,
symmetric functions and Grothendieck polynomials}, arXiv:hep-th/9306005.

\bibitem[Ga]{gasharov} {V. Gasharov}, {\it A short proof of the
Littlewood-Richardson rule}, {European J. Combin.} {\bf 19} (1998) 451--453.

\bibitem[JP]{jp:ss} {G. D. James and M. H. Peel}, {\it Specht series for
skew representations of symmetric groups}, {J. Algebra} {\bf 56} (1979),
343--364.

\bibitem[LR]{lr:gc} {D. E. Littlewood and A. R. Richardson},
{\it Group characters and algebra,}
{Philos. Trans. Roy. Soc. London Ser. A} {\bf 233} (1934), 49--141.


\bibitem[La]{ldouble}
{A. Lascoux}, {``Interpolation",\/} Lectures at Tianjin
University, June 1996.

\bibitem[LS]{lassch}
{A. Lascoux and M. P. Sch\"utzenberger}, {\em Structure de Hopf de
l'anneau de cohomologie et de l'anneau de Grothendieck d'une
vari\'et\'e de drapeaux}, C.R. Acad. Sci. Parix S\'er. I Math,
{\bf 295} (1982), 629--633.

\bibitem[Le]{Lenart} {C. Lenart}, {\it Combinatorial Aspects of the
$K$-Theory of Grassmannians}, Ann. Comb, {\bf 4} (2000), 67--82.

\bibitem[Mac1]{Mc:Schur} {I. G. Macdonald}, {\it Schur functions: theme and
variations}, { in} {``Actes 28-e S\'eminaire Lotharingien",\/}
{Publ. I.R.M.A. Strasbourg,\/} 1992, 498/S--27, 5--39.

\bibitem[Mac2]{Mc}
{I. G. Macdonald}, {\em ``Symmetric Functions and Hall
Polynomials,"} 2nd edition, Oxford University Press, Oxford 1995.

\bibitem[MS]{MS:main}
{A. I. Molev and B. E. Sagan}, {\em A Littlewood-Richardson Rule
For Factorial Schur Functions}, Trans. Amer. Math. Soc, {\bf 351}
(1999), 4429--4443.

\bibitem[Mo]{M:vanis}
A. Molev, \textit{Factorial Supersymmetric Schur Functions and
Super Capelli Identities}, Amer. Math. Soc. Transl, {\bf 181}
(1998), 109--137.

\bibitem[Ok]{M:vanish} A. Okounkov, {\it Quantum immanants and higher
Capelli identities}, {Transformation Groups} {\bf 1} (1996)
99--126.

\bibitem[Sa]{sagan} {B. E. Sagan}, {\em ``The symmetric group:
representations, combinatorial algorithms, and symmetric
functions,"} 2nd edition, Springer, New York 2001.

\bibitem[Ze]{z:gl} {A. V. Zelevinsky}, {\it A generalization
of the Littlewood--Richardson rule and the Robinson--Schensted--Knuth
correspondence},
 {J. Algebra} {\bf 69} (1981),
82--94.

\end{thebibliography}

\end{document}